\documentclass{article}

\RequirePackage{amsthm,amsmath,amsfonts,amssymb}
\RequirePackage{mathrsfs,bbm} 
\RequirePackage{natbib}
\RequirePackage{authblk}
\RequirePackage{hyperref}
\RequirePackage[margin=26mm]{geometry}

\theoremstyle{plain}

\newtheorem{theorem}{Theorem}[section]
\newtheorem{lemma}[theorem]{Lemma}
\newtheorem{prop}[theorem]{Proposition}
\newtheorem{coro}[theorem]{Corollary}
\theoremstyle{definition}
\newtheorem{definition}[theorem]{Definition}
\newtheorem*{example}{Example}

\newtheorem{assumption}[theorem]{Assumption}
\theoremstyle{remark}

\DeclareMathOperator*{\argmax}{arg\,max}
\DeclareMathOperator*{\argmin}{arg\,min}

\def\namedlabel#1#2{\begingroup
    #2%
    \def\@currentlabel{#2}%
    \phantomsection\label{#1}\endgroup
}

\title{Non-asymptotic analysis of label shift quantification: robustness and calibration\thanks{This is a revised version of the manuscript submitted to TEST.}}
\author[1]{Alexandre Lecestre}
\affil[1]{\footnotesize{University of Bordeaux, CNRS, Bordeaux INP, IMB, UMR 5251, F-33400 Talence, France}}

\renewcommand\footnotemark{}

\begin{document}

\maketitle

\begin{abstract}
In this paper, we investigate the label shift quantification problem. We propose robust estimators of the label distribution which turn out to coincide with the Maximum Likelihood Estimator. We analyze the theoretical aspects and derive deviation bounds for the proposed method, providing optimal guarantees in the well-specified case, along with notable robustness properties against outliers and contamination. Our results provide theoretical validation for empirical observations on the robustness of Maximum Likelihood Label Shift.
\end{abstract}


\section{Introduction}\label{sec1}

The assumption that training and test samples share the same data generation process is at the base of most supervised learning methods. However, this assumption often fails in real-world applications, posing challenges to practitioners. The following example inspired by the introduction of \cite{kelipton_wang_smola} is a good illustration of the problem.
\begin{example}
We have been studying a specific disease (say cholera/hepatitis A-E) and training classifier $\eta$ to detect whether a person is suffering from the disease based on well chosen covariates $x$, where $\eta(x)=1$ predicts a diseased patient and $\eta(x)=0$ predicts a healthy patient. This was performed under 'normal' conditions, where the proportion of diseased individuals in the training set is $\alpha>0$. During an epidemic, the proportion of diseased individuals being tested becomes significantly higher. This violates the common \emph{i.i.d.} (independent and identically distributed observations) assumption and will render the classifier inefficient as it will underestimate the diseased rate. This underestimation arises because the classifier, trained on a lower prevalence of the disease, assumes the same proportions hold in the test data.
\end{example}
This example illustrates a common real-world challenge where the training and test datasets do not follow the same distribution—a phenomenon known as distribution shift. To address this, the classifier must be adapted to the new data. However, achieving this is infeasible without specific assumptions about the nature of the shift. This paper focuses on label shift, commonly used in classification contexts, assuming the conditional distribution of the covariates remains unchanged between the training and test datasets. In this scenario the training dataset is labeled but the test dataset is not. It differs from covariate shift, naturally used in prediction or regression contexts where covariates $x$ cause a response $y$, which assumes that the conditional of the response variable is the same in both the training and test samples. Let us describe more formally the label shift assumption.\\
Let $\mathscr{X}$ denote the covariate space and $\mathscr{Y}:=\{1,2,\dots,k\}=:[k]$ the label space, where $k$ is an integer larger than 1. We denote by $D_s(dx,dy)$, respectively $D_t(dx,dy)$, the distribution of the training data over $\mathscr{X}\times\mathscr{Y}$, respectively the test data. In the literature, $D_s$ and $D_t$ are sometimes called \emph{source domain} and \emph{target domain}. The label shift assumption corresponds to
\[ D_s(dx|y) = D_t(dx|y) \text{  for all } y\in\mathscr{Y}. \]
In the example above, it means that the symptoms of the disease have not changed between the training period and now, only the proportions have changed. This assumption enables the adaptation of our predictor to new data, without the need to train a new model from scratch.\par
The source label distribution $\alpha^*\in\mathcal{W}_k$ and the target label distribution $\beta^*\in\mathcal{W}_k$ are given by
\[ \alpha^*:= (D_s(\mathscr{X}\times\{i\}))_{i\in[k]} \text{  and  } \beta^*:=(D_t(\mathscr{X}\times \{i\}))_{i\in[k]}, \]
where
\[ \mathcal{W}_k = \{ x\in[0,1]^k; x_1+\dots+x_k =1 \} \]
is the simplex. The simplex $\mathcal{W}_k$ is identified as the class of all probability distributions over $[k]$ here and in the rest of the paper. The literature addresses several challenges within the context of label shift. Detection involves determining if there has been a distribution shift, i.e. testing whether $\beta^*=\alpha^*$ or $\beta^*\neq\alpha^*$. Correction aims to produce a classifier that performs well for the target distribution $D_t$. In this paper, we focus on a third problem called label shift quantification, or label shift estimation, where the objective is to estimate the target label distribution $\beta^*$, or similarly, the vector of ratios $w^*=(\beta^*_i/\alpha^*_i)_{i\in[k]}$. There is a rich body of literature on the subject. We refer the reader to \cite{dussap,garg,alexandari} for extensive introductions to the topic.\\
A naive method is to estimate the conditional distributions on the training data and the target label distribution on the test data. It can be done using common estimators developed for mixture models as there is a natural relation between label shift and mixtures. A finite mixture distribution is a distribution of the form
\[ w_1 F_1 + \dots + w_K F_K, \] where $K\geq 2$ is an integer, $w=(w_1,\dots,w_K)\in\mathcal{W}_K$ is the vector of weights and $F_1,\dots,F_K$ are probabilities called emission distributions. Under the label shift assumption, the source and target covariate distributions can be expressed as mixtures, i.e.
\[ D_s(dx,[k]) =  \alpha^*_1 Q^*_1 + \dots + \alpha^*_k Q^*_k \text{ and  } D_t(dx,[k]) = \beta^*_1 Q^*_1 + \dots + \beta^*_k Q^*_k, \]
where $Q ^*_i(dx)=D_s(dx|i)=D_t(dx|i)$ for all $i\in[k]$. Each distribution $Q ^*_i$ represents the distribution of covariates conditioned on label $i$, and the label proportions determine the overall distribution. Both are $k$-component mixture distributions sharing with the same emission distributions $Q^*_1,\dots,Q^*_k$. For the label shift quantification problem to be well-posed, we need a linear independence assumption on the conditional distributions $Q^*_1,\dots,Q^*_k$. Otherwise, the problem we are considering is ill-posed as the vector $\beta^*$ is not identifiable. The strategy we mentioned earlier, i.e. to estimate $Q^*_1,\dots,Q^*_k$ on the source dataset and later $\beta^*$ on the target dataset, works well in theory, but it requires rather accurate estimates of $Q^*_1,\dots,Q^*_k$ which is challenging in practice. This approach becomes impractical in high-dimensional settings, where the covariate space $\mathscr{X}$ is large relative to the sample size.\\
We can mention some common methods that have been developed to address this problem, such as distribution matching using Reproducing Kernel Hilbert Spaces which was inspired by the Kernel Mean Matching (KMM) approach (see \cite{pmlr-v32-iyer14}). Another one is Black Box Shift Estimation (BBSE), where one uses the confusion matrix of an off-the-shelf classifier to adjust the predicted label distribution (see \cite{kelipton_wang_smola}). Maximum Likelihood Label Shift (MLLS), probably the most common approach, applies the maximum likelihood principle using an off-the-shelf classifier. The different methods have been widely studied through theoretical guarantees and empirical performances. Although they appear as distinct approaches, recent papers seem to reveal similarities. \cite{garg} established the theoretical equivalence of the optimization objectives in MLLS and BBSE. Similarly, \cite{dussap} introduce Distribution Feature Matching (DFM) which is a general framework including KMM and BBSE. \cite{dussap} also extend the classical framework to consider the \emph{contaminated label shift} setting where the covariate distribution of the test sample is of the form
\begin{equation}
\label{eq:contaminated_label_shift}
\beta^*_0 Q_0 + \beta^*_1 Q^*_1 + \dots + \beta^*_k Q^*_k,
\end{equation}
with $\beta^*\in\mathcal{W}_{k+1}$, modelling a contaminated dataset under label shift. In this setting, they are interested in estimating the weights $(\beta^*_i)_{i\geq 0}$ and obtain a general result with Corollary 1 in \cite{dussap}. Their deviation bound shows that their estimator is robust, but only to a specific type of contamination. To our knowledge, it is the only theoretical guarantees of robustness in the label shift settings and we aim to fill this gap.\\
In this work, we propose robust methods for the estimation of the target label distribution $\beta^*$. It turns out that our method includes maximum likelihood approaches such that our results apply, in particular, to the MLE. We consider two different scenarios, where we build our strategy upon off-the-shelf estimators in both cases. In the first one, we are given estimates of the conditional distributions. It is related to weight vector estimation, as studied in \cite{dalalyan_ems} and \cite{tsybakov_spades}, with notable differences in the approach. The second scenario corresponds to the setting of label shift, as described in \cite{garg}, and in this case, we have a predictor that was trained on the source domain.\\
We provide a thorough theoretical analysis of our estimation strategies, including general deviation bounds under minimal assumptions and establish convergence rates in well-specified settings. Furthermore, we investigate the robustness of our estimators to misspecification, contamination, and outliers. The contamination setting described above--see (\ref{eq:contaminated_label_shift})--corresponds to the standard Huber contamination model. We leave the study of the stronger adversarial contamination setting (REF) for future work. Note that we consider the \emph{contaminated label shift setting} (\ref{eq:contaminated_label_shift}) but with a goal different from \cite{dussap}. We do not aim to estimate $(\beta^*_i)_{i\geq 0}$ in general but we want our estimator of $(\beta^*_i)_{i\geq 1}$ to be robust to small deviations from the ideal framework, i.e. small values of $\beta^*_0$ in this case. Indeed, we show that our estimator's performance depends solely on the contamination rate, regardless of its nature. In practice, datasets are often noisy or contain outliers, making robustness a crucial property for reliable estimation. These results of robustness complete previous works on MLLS. We are bridging further the gap between prior empirical studies (e.g. \cite{saerens,alexandari}) and the theoretical results of \cite{garg}.
\section{Statistical framework}
Let $\mathscr{X}$ be the covariate space, endowed with a $\sigma$-algebra $\mathcal{X}$, such that $\left(\mathscr{X},\mathcal{X}\right)$ is a measurable space. We denote by $\mathscr{P}_X$ the class of all probability distributions on $\left(\mathscr{X},\mathcal{X}\right)$. Let $X_1,\dots,X_n$ be independent random variables on $\left(\mathscr{X},\mathcal{X}\right)$. Those random variables correspond to the covariates in the target data, where $n$ is the size of the target sample size. We denote by $P_i\in\mathscr{P}_X$ the distribution of the random variable $X_i$ for all $i\in[n]$. We will often work under the following assumption.
\begin{assumption}
\label{hyp:iid_mix}
The variables $X_1,\dots,X_n$ are \emph{i.i.d.} with common distribution $P^*$ of the form
\begin{equation}
\label{eq:assumption_iid_mix}
P^* := \beta^*_1 Q^*_1 + \dots + \beta^*_k Q^*_k, 
\end{equation}
where $k\geq 2$ and $\beta^*\in\mathcal{W}_k$. Moreover, $Q^*_1,\dots,Q^*_k$ are linearly independent in the space of signed measures on $\left(\mathscr{X},\mathcal{X}\right)$.
\end{assumption}
This assumption means that the observations are i.i.d. with a common distribution which can be written as a finite mixture. The integer $k$ is known and corresponds to the number of different labels observed in the training phase. We always consider that the sample size satisfies $n\geq e\times k$. The linear independence of $Q^*_1,\dots,Q^*_k$ ensures that the label distribution $\beta^*\in\mathcal{W}_k$ is identifiable. It is further discussed in the first part of Section \ref{sec:mpe}.\par
However, we do not want to rely on this assumption for our estimation strategy to be effective. In practice, the data may not be perfectly i.i.d., and the distribution of the observations may deviate from a finite mixture model. We will obtain general results when we do not make any assumption on those distributions $P_1,\dots,P_n$. This will allow us to consider the possible presence of contamination or outliers and quantify the robustness of our estimator in these cases. Our estimation strategy is to do as if Assumption \ref{hyp:iid_mix} were true, but is designed to remain effective even when this assumption is violated.\par
We consider two different settings. The second setting is more specific to the label shift problem, but the first one is more direct and allows us to introduce and discuss notions necessary to consider the second setting.
\begin{itemize}
\item[Setting \namedlabel{case:1}{A}] We already know the conditional distributions $(Q^*_i)_{i\in[k]}$ from the source dataset. If not, we are given estimates $(Q_i)_{i\in[k]}$ of those distributions. We will always assume that the distributions $Q_1,\dots,Q_k$ are linearly independent.
\item[Setting \namedlabel{case:predictor}{B}] We know the Bayes predictor $f^*:\mathscr{X}\rightarrow \mathcal{W}_k$ -- explicitly given in (\ref{eq:f_bayes}) -- for the source distribution. If not, we are given estimates $f$ of $f^*$ and  $\alpha$ of $\alpha^*$.
\end{itemize}
We do not consider the problem of estimating the quantities $Q^*_1,\dots,Q^*_k$ and $f^*$ here as they have already been widely investigated in the literature. We call predictor any measurable function $\mathscr{X}\rightarrow\mathcal{W}_k$, where we assume that $\mathcal{W}_k$ is naturally endowed with the $\sigma$-algebra induced by the Borel $\sigma$-algebra on $\mathbb{R}^d$. We prefer to work with a predictor, giving label probabilities, rather than with a hard classifier, i.e. a function $\mathscr{X}\rightarrow[k]$. A classifier $g$ can always be deduced from a predictor $f$ using the label with maximum probability, i.e. $g(x)=\argmax_{i\in[k]} f_i(x)$. However, the predictor carries more information than a classifier which is crucial to perform label shift estimation. In our context, investigating the first setting is more direct and allows us to introduce and discuss notions necessary to consider the second setting. Therefore, Sections \ref{sec:mpe} and \ref{sec:main} correspond to Settings \ref{case:1} and \ref{case:predictor} respectively.\\
Our estimation strategy is based on $\rho$-estimators introduced by \cite{rho_inventiones,baraudrevisited}. It is a model-based estimation method which is proven to be robust to small deviations, those deviations being quantified via the Hellinger distance. The Hellinger distance between two distributions $Q$ and $Q'$ on the same measurable space is defined by
\[ h^2(Q,Q') = \frac{1}{2} \int \left(\sqrt{\mathrm{d}Q/\mathrm{d}\nu} -\sqrt{\mathrm{d}Q'/\mathrm{d}\nu}\right)^2 \mathrm{d}\nu, \]
where $\nu$ is any positive measure that dominates both $Q$ and $Q'$, the result being independent of $\nu$. The Hellinger distance is particularly appealing from a robustness perspective, especially compared to the Kullback-Leibler (KL) divergence, which is intrinsically linked to the maximum likelihood approach. The KL divergence is finite only if the true distribution is absolutely continuous with respect to distributions in our model. This can be problematic in the presence of contamination by an atypical distribution. In contrast, the Hellinger distance is always well-defined, remains bounded by 1, and has the additional advantage of being symmetric.\\
Next, we define $\rho$-estimators, which are naturally suited for addressing our estimation problem. For a quick reading of this paper, the reader may skip to the next section and simply think of our estimator as the MLE.
\subsection{\texorpdfstring{$\rho$}{rho}-estimation}
\label{sec:rho_est}
Let $\mathscr{Q}$ be a countable subset of $\mathscr{P}_X$, and $\mathcal{Q}$ be the associated set of densities with respect to a $\sigma$-finite product measure $\mu$ -- see \cite[page 9]{legall} -- on $\left(\mathscr{X},\mathcal{X}\right)$, such that
\[ \mathscr{Q} = \{ q\cdot \mathrm{d}\mu: q\in\mathcal{Q} \}, \]
where $q\cdot \mathrm{d}\mu$ denotes the measure with density $q$ with respect to $\mu$. We define $\rho$-estimators on $\mathscr{Q}$ as follows. We denote by $\psi$ the function defined by 
\begin{equation}
\label{eq:psi_rho}
\psi: \begin{array}{|lcl}
[0,+\infty] & \rightarrow & [-1,1]\\
x & \mapsto & \frac{x-1}{x+1}
\end{array},
\end{equation}
with the convention $\psi(+\infty)=1$. For $\mathbf{x} = (x_1,...,x_n) \in \mathscr{X}^n$ and $q,q'\in \mathcal{Q}$, we define
\[ \mathbf{T}(\mathbf{x},q,q') := \sum\limits_{i=1}^{n}\psi\left(\sqrt{\frac{q'\left(x_i\right)}{q\left(x_i\right)}}\right), \]
with the convention $0/0 = 1$ and $a/0 = +\infty$ for all $a >0$. To build an intuition on why are $\rho$-estimators robust and why this definition, one should see $\mathbf{T}$ as a robust version of the likelihood ratio test (LRT). For instance, if you were to take $\psi=\log$, you would fall back on the LRT. We refer to Proposition 2 and (12) in \cite{baraudrevisited} for guarantees on the test $\mathbf{T}$. We define
\[ \mathbf{\Upsilon}(\mathbf{x},q) := \sup_{q'\in\mathcal{Q}} \mathbf{T}(\mathbf{x},q,q'), \]
for all density $q\in\mathcal{Q}$. For random variables $X_1,\dots,X_n$ on $(\mathscr{X},\mathcal{X})$, the $\rho$-estimator $\hat{P}(\mathbf{X},\mathcal{Q})$ is any measurable element of the closure (with respect to the Hellinger distance) of the set
\begin{equation}
\label{eq:rho_est}
\pmb{\mathscr{E}}(\mathbf{X}) := \left\{ Q = q\cdot \mathrm{d}\mu ; q\in\mathcal{Q}, \mathbf{\Upsilon}(\mathbf{X},q) < \inf_{q'\in\mathcal{Q}} \mathbf{\Upsilon}(\mathbf{X},q')  + 11.36 \right\}.
\end{equation}
The constant $11.36$ is given by (7) in \cite{baraudrevisited}. It does not play an essential role and can be replaced by any smaller positive constant, at the cost of a larger constant in the theoretical guarantees. It is a technical artefact to ensure $\rho$-estimators are well defined. We require $\mathscr{Q}$ to be countable for the same reason, but in practice it will just mean that we consider only rational parameters in $\mathbb{Q}^d$ instead of all the real valued parameters in $\mathbb{R}^d$. Since we take the closure with respect to the Hellinger distance in the definition, $\rho$-estimators can also correspond to parameters in $\mathbb{R}^d\backslash\mathbb{Q}^d$.
\section{Label shift quantification in Setting \ref{case:1}}
\label{sec:mpe}
This setting has attracted considerably less attention than Setting \ref{case:predictor} in the label shift literature, though it is closely connected to the problem of estimating mixing weights in mixture models. Relevant works include those of Dalalyan \& Sebbar \cite{dalalyan_ems} and Bunea et al. \cite{tsybakov_spades}, who study the estimation of a target density $f^*$ using a mixture $f_{\lambda}=\sum\limits_{i=1}^M \lambda_i f_i$. Here, $(f_i)_{1\leq i\leq M}$ is a fixed dictionary of densities, and the setting assumes sparsity $(M\gg n)$, where $f^*$ can be well-approximated by $f_{\lambda}$ with most $\lambda_i$ vanishing. The cited works establish theoretical guarantees not only for the estimation of $f^*$ but also for the recovery of the mixing weights $\lambda$ giving the best approximation of $f^*$.
\subsection{Our estimator}
We introduce here our estimation strategy for the first setting (\ref{case:1}). Let $Q_1,\dots,Q_k$ be (linearly independent) distributions in $\mathscr{P}_X$ and $q_1,\dots,q_k$ be their associated densities with respect to a $\sigma$-finite measure $\mu$. We define the mixture model 
\[ \mathscr{M}_{mix}(Q_1,\dots,Q_k) := \left\{ \sum\limits_{i=1}^k \beta_i Q_i; \beta \in\mathcal{W}_k\cap\mathbb{Q}^k \right\}, \]
which is a countable and dense subset of
\[ \overline{\mathscr{M}}_{mix}(Q_1,\dots,Q_k) := \left\{ \beta_1 Q_1 + \dots + \beta_k Q_k; \beta \in\mathcal{W}_k \right\}, \]
with respect to the Hellinger distance. We denote the associated class of densities by
\begin{equation}
\label{eq:mix_fix_dens}
\mathcal{M}_{mix}(q_1,\dots,q_k) := \left\{ \sum\limits_{i=1}^k \beta_i q_i; \beta \in\mathcal{W}_k\cap\mathbb{Q}^k \right\}.
\end{equation}
Before presenting the results, we discuss our estimation strategy. Let us put ourselves in the context of Assumption \ref{hyp:iid_mix}. Our method is to build an estimator $\hat{\beta}$ of $\beta^*$ from a $\rho$-estimator $\hat{P}=\hat{P}(\mathbf{X},\mathcal{M}_{mix}(q_1,\dots,q_k))$, as defined by (\ref{eq:rho_est}). The next result establishes the connection with the maximum likelihood approach in this context.
\begin{prop}
\label{prop:mle}
When it exists, the Maximum Likelihood Estimator $\hat{\beta}_{MLE}$ given by
\begin{equation}
\label{eq:mle_a}
\hat{\beta}_{MLE} \in \argmax_{\beta\in\mathcal{W}_k} \sum\limits_{i=1}^n \log\left( \sum\limits_{j=1}^k \beta_j q_i(X_j) \right)
\end{equation}
is a $\rho$-estimator with respect to $\mathcal{M}_{mix}(q_1,\dots,q_k)$.
\end{prop}
This result is proven in Section \ref{sec:proof_prop_mle}. It implies that all the results we will give for our estimator are also valid for the maximum likelihood estimator $\hat{\beta}_{MLE}$. We are only aware of the article of Dalalyan \& Sebbar \cite{dalalyan_ems} which considers the MLE in a similar setting, but with no relation to label shift and no robustness considerations. One difference is that the MLE might not exist while $\rho$-estimators are always well-defined. For instance, we do not need to assume that the considered densities are bounded. Another implication of this result is that standard methods like the EM-algorithm can be used to compute our estimator.\\
In the case of $\rho$-estimation, we still need to give a way to deduce our estimator $\hat{\beta}$ from a $\rho$-estimator $\hat{P}=\hat{P}(\mathbf{X},\mathcal{M}_{mix}(q_1,\dots,q_k))$. We also need to make sure it is the right approach. If $Q_i$ is relatively close $Q^*_i$ for all $i$, the model $\overline{\mathscr{M}}_{mix}(Q_1,\dots,Q_k)$ is a good approximation of $P^*$. Indeed, we have
\[ h\left(P^*,\beta^*_1 Q_1 + \dots + \beta^*_k Q_k\right) \leq \max_{i\in[k]} h(Q^*_i,Q_i), \]
and this is due to the following result.
\begin{lemma}{(Lemma B.3 \cite{lecestre_mixture})}\\
For all $w,w'\in\mathcal{W}_k$ and all $F_1,F'_1,\dots,F_k,F'_k$ in $\mathscr{P}_X$, we have
\begin{equation}
\label{eq:mix_hellinger_upper}
h\left( \sum\limits_{i=1}^k w_i F_i, \sum\limits_{i=1}^k w'_i F'_i \right) \leq h(w,w') + \max_{i\in[k]} h(F_i,F'_i).
\end{equation}
\end{lemma}
This means we can obtain a good estimator of $P^*$ as long as $\max_{i\in[k]} h(Q^*_i,Q_i)$ is small (see Lemma \ref{lem:aux_2}). However, our goal is to estimate the target label distribution $\beta^*$. One would naturally consider any $\hat{\beta}\in\mathcal{W}_k$ satisfying $\hat{P}=\hat{\beta}_1 Q_1 + \dots + \hat{\beta}_k Q_k$ as an estimator of $\beta^*$. One issue is that the vector $\hat{\beta}$ is not uniquely defined when the distributions $Q_1,\dots,Q_k$ are not linearly independent. This is why we assume the distributions $Q_1,\dots,Q_k$ to be linearly independent. This assumption ensures the identifiability of the mixture weights,
that is, the representation (8) is unique.
\\
From now on, our estimator $\hat{\beta}$ of $\beta^*$ is defined as the unique element of $\mathcal{W}_k$ such that
\begin{equation}
\label{eq:hat_beta}
\hat{P} = \hat{\beta}_1 Q_1 + \dots + \hat{\beta}_k Q_k,
\end{equation}
where $\hat{P}=\hat{P}(\mathbf{X},\mathcal{M}_{mix}(q_1,\dots,q_k))$ is a $\rho$-estimator. 
Inequality (\ref{eq:mix_hellinger_upper}) indicates that if $\hat{\beta}$ is close to $\beta^*$ then $\hat{P}$ is close to $P^*$ but we need to obtain the converse. The next result shows that it is true under the linear independence assumption.
\begin{lemma}
\label{lem:aux_1}
For all distributions $F_1,\dots,F_k$ in $\mathscr{P}_X$ and all $\beta,\overline{\beta}\in\mathcal{W}_k$ we have
\[ h\left( \sum\limits_{i=1}^k \beta_i F_i, \sum\limits_{i=1}^k \overline{\beta}_i F_i \right) \geq \frac{\Delta^*(F_1,\dots,F_k)}{2\sqrt{2}} ||\beta-\overline{\beta}||_1, \]
where
\[ \Delta^*(F_1,\dots,F_k) = \inf_{\substack{I\subset [k]\\I\neq[k]}} \inf_{\gamma\in\mathcal{W}_{|I|}} \inf_{\lambda\in\mathcal{W}_{k-|I|}} d_{TV}\left( \sum\limits_{i\in I} \gamma_i F_i, \sum\limits_{i\in [k]\backslash I} \lambda_i F_i \right), \]
where $d_{TV}$ is the total variation distance and $|I|$ denotes the cardinal of the set $I$.
\end{lemma}
The proof can be found in Section \ref{sec:proof_lem_aux1}. One can check that $\Delta^*(F_1,\dots,F_k)$ is a positive constant as soon as the distributions $F_1,\dots,F_k$ are linearly independent. The quantity $\Delta^*(F_1,\dots,F_k)$ measures how well separated the components of the mixture are. It is possible to compute this constant from $F_1,\dots,F_k$ but it should be easier to compute a lower bound on $\Delta^*(F_1,\dots,F_k)$ if we have associated densities $f_1,\dots,f_k$ (with respect to a $\sigma$-finite measure $\mu$) that are bounded, e.g. by a constant $M$. In that case, we have
\[ \Delta^*(F_1,\dots,F_k) \geq \frac{1}{2M} \inf_{\substack{I\subset [k]\\I\neq [k]}} \inf_{\gamma\in\mathcal{W}_{|I|}} \inf_{\lambda\in\mathcal{W}_{k-|I|}} \left|\left|\sum\limits_{i\in I} \gamma_i f_i - \sum\limits_{i\in [k]\backslash I}\lambda_i f_i \right|\right|_{L_2(\mu)}^2, \]
and finding the right hand side of this inequality is a quadratic programming problem.\par
Results of label shift quantification can take different forms. Some of them consider the $\ell_2$-distance from the $\hat{\beta}$ to the target $\beta^*$ \cite{dussap}, other consider the $\ell_2$-distance from the vector $w$ of shift ratios, given by $w_i=\beta^*_i/\alpha^*_i$, to its estimate counterpart \cite{garg}. In this paper, we consider the $\ell_1$-distance from $\hat{\beta}$ to $\beta^*$, which is a natural choice from a probabilistic point of view as it corresponds to the total variation distance between distributions in $\mathcal{W}_k$. This can make the comparison of different results a bit more difficult, but we can still mention the inequalities 
\begin{align*}
||\beta-\beta'||_2 \leq &||\beta-\beta'||_1 \leq \sqrt{k} ||\beta-\beta'||_2\\
\left(\min\limits_{1\leq i\leq k} \alpha_i\right) ||w-w'||_2 \leq &||\beta-\beta'||_1 \leq ||\alpha||_1 ||w-w'||_2,
\end{align*}
which hold for all $\beta,\beta'\in\mathcal{W}_k$, where $w_i=\beta_i/\alpha_i$ and $w'_i=\beta'_i/\alpha_i$ for all $i\in[k]$.
\subsection{Results}
Following Lemma \ref{lem:aux_1}, we can obtain a deviation inequality for our estimator $\hat{\beta}$ since the distributions $Q_1,\dots,Q_k$ are assumed to be linearly independent. In the rest of this section, $\hat{\beta}$ will denote either the estimator given by (\ref{eq:hat_beta}) or the MLE defined by (\ref{eq:mle_a}).
\begin{theorem}
\label{th:conditional}
\begin{itemize}
\item There is a positive constant $C(\mathbf{Q})$ depending only on\\
$Q_1,\dots,Q_k$, such that for all $\overline{\beta}\in\mathcal{W}_k$ and all $\xi>0$,
\begin{equation}
\label{eq:th_conditional}
C(\mathbf{Q}) ||\overline{\beta}-\hat{\beta}||_1^2 \leq n^{-1} \sum\limits_{j=1}^n h^2\left( P_j, \sum\limits_{i=1}^k \overline{\beta}_i Q_i  \right) + \frac{k \log (n/k) + \xi}{n},
\end{equation}
with probability at least $1-e^{-\xi}$.
\item Let $\overline{Q}_1,\dots,\overline{Q}_k$ be linearly independent distributions in $\mathscr{P}_X$. There is a positive constant $C(\overline{\mathbf{Q}})$ depending only on $\overline{Q}_1,\dots,\overline{Q}_k$ such that for all $\overline{\beta}\in\mathcal{W}_k$ and all $\xi>0$,
\begin{align*}
C(\overline{\mathbf{Q}}) ||\overline{\beta}-\hat{\beta}||_1^2 &\leq n^{-1} \sum\limits_{j=1}^n h^2\left( P_j, \sum\limits_{i=1}^k \overline{\beta}_i \overline{Q}_i \right) + \max_{1\leq i\leq k} h^2\left(Q_i,\overline{Q}_i\right)\\
&+ \frac{k \log (n/k) + \xi}{n},
\end{align*}
with probability at least $1-e^{-\xi}$.
\end{itemize}
\end{theorem}
This result is proven in Section \ref{sec:proof_th_conditional}. This result holds without Assumption \ref{hyp:iid_mix}. Therefore, we need to consider any $\overline{\beta}$ in $\mathcal{W}_k$ as $\beta^*$ does not necessarily exist. Corollary \ref{coro:conditional} below is the result we obtain under Assumption \ref{hyp:iid_mix}.\par
The second inequality of Theorem \ref{th:conditional} allows us to have a constant depending depending on other conditional distributions. Typically, if $Q_1,\dots,Q_k$ are estimators of $Q^*_1,\dots,Q^*_k$ under Assumption, the constant $C(Q)$ in the first inequality is random as it depends on the training dataset used to train $Q_1,\dots,Q_k$. This issue is avoided with the second inequality. Theorem \ref{th:conditional} is a very general result and is not very informative without any assumption on the distributions $P_1,\dots,P_n$. The quantity 
\[ \sum\limits_{j=1}^n h^2\left( P_j, \sum\limits_{i=1}^k \overline{\beta}_i Q_i  \right) \]
in (\ref{eq:th_conditional}) quantifies the distance from our model to the true distributions of the observations. It needs to be sufficiently small for our bound to be meaningful but not equal to 0 necessarily. It happens in the \emph{contaminated label shift} setting described by (\ref{eq:contaminated_label_shift}) with a low contamination rate $\beta^*_0$ for example. Before discussing the robustness properties of the estimator, we illustrate Theorem \ref{th:conditional} with the following result, which is a direct consequence of the second deviation inequality therein.
\begin{coro}
\label{coro:conditional}
Under Assumption \ref{hyp:iid_mix}, there is a positive constant $C(\mathbf{Q}^*)$ depending only on $Q^*_1,\dots,Q^*_k$ such that for all $\xi>0$,
\[ C(\mathbf{Q}^*) ||\beta^*-\hat{\beta}||_1^2 \leq \max_{1\leq i\leq k} h^2\left( Q_i, Q^*_i \right) + \frac{k \log (n/k) + \xi}{n}, \]
with probability at least $1-e^{-\xi}$. It is possible to replace the constant $C(\mathbf{Q}^*)$ by a constant $C(\mathbf{Q})$ depending on $Q_1,\dots,Q_k$ instead.
\end{coro}
This result offers key insights into the convergence of the mixture proportion estimator $\hat{\beta}$. Under mild assumptions on the linear independence of the component distributions, Corollary \ref{coro:conditional} ensures that $\hat{\beta}$ achieves a near-parametric convergence rate of $n^{-1/2} \log^{1/2} n$ under well-specified settings, with respect to the $\ell_1$-loss. This rate is comparable (up to logarithmic factors) to that of \cite[Proposition 3.2]{dalalyan_ems} which also analyzes maximum likelihood estimation. While our rate is slightly slower due to the $\log n$ term, our assumptions are weaker--notably, we avoid restrictive conditions like equation (1.12) in \cite{dalalyan_ems}. In contrast, the bound of Bunea \emph{et al.} \cite{tsybakov_spades} is not suited to label shift, as their framework targets sparse regimes with a very large dictionary $M\gg n$.\\
Though these works establish guarantees the mixing weights, in particular for the MLE in \cite{dalalyan_ems}, they do not address robustness--a central focus of our analysis. We now consider this aspect for our estimator with specific cases of interest. The case of misspecification actually corresponds to Corollary \ref{coro:conditional}. One can notice that the performance of our estimator does not degrade notably as long as the misspecification term $\max_{i\in[k]}h^2(Q_i,Q^*_i)$ is of order not larger than $n^{-1} \log n$.\\
The Huber contamination model corresponds to the case studied by \cite{dussap} in which $X_1,\dots,X_n$ are \emph{i.i.d.} with common distribution $P^*$ given by
\begin{equation}
\label{eq:p_contamination}
P^* = \lambda_0 \overline{P} + \lambda_1 Q^*_1 + \dots + \lambda_k Q^*_k,
\end{equation}
where $\overline{P}$ is any distribution in $\mathscr{P}_X$. Our method is designed to retrieve the 'original weights' $\beta^*$ we had before contamination given by $\beta^*_i=\lambda_i/(1-\lambda_0)$, for all $i\in[k]$. Assuming $Q_i=Q^*_i$ for all $i\in[k]$, we have
\begin{align*}
C(\mathbf{Q}^*) ||\beta^*-\hat{\beta}||_1^2 &\leq \lambda_0 + \frac{ k \log (n/k) + \xi}{n},
\end{align*}
with probability at least $1-e^{-\xi}$. As long as $\lambda_0$ is small compared to $n^{-1} k \log(n/k)$, the performance of our estimator is not significantly worse than in the ideal case without contamination.\\
We model the presence of outliers in the following way. Assume $X_1,\dots,X_n$ are independent and there is an index set $I\subset [n]$ of outliers, i.e.
\begin{equation}
\label{eq:p_outliers}
X_i\sim P^*=\beta^*_1 Q^*_1 + \dots + \beta^*_k Q^*_k \text{  for all  } i\in[n]\backslash I
\end{equation}
and $X_i$ follows any distribution $P_i$ in $\mathscr{P}_X$ for $i\in I$. In that case, assuming $Q_i=Q^*_i$ for all $i\in[k]$, we have
\begin{align*}
C(\mathbf{Q}) ||\beta^*-\hat{\beta}||^2_1 &\leq \frac{|I|}{n} + \frac{k \log (n/k) + \xi}{n},
\end{align*}
with probability at least $1-e^{-\xi}$, for all $\xi>0$. As long as the proportion of outliers $|I|/n$ is small compared to $n^{-1} k \log(n/k)$, the performance of the estimator is still of the same order as in the ideal case without contamination.\\
The approach developed in this section relies on obtaining good estimates of the conditional distributions $Q^*_1,\dots,Q^*_k$. However, this becomes increasingly difficult when these distributions belong to high-dimensional models and we cannot hope for a reasonably small value of $\max_{1\leq i\leq k} h^2(\hat{Q}_i,Q^*_i)$ in the bound of Corollary \ref{coro:conditional}. In addition, while the study of the distributions $Q^*_1,\dots,Q^*_k$ may be of interest, depending on the context, their estimation is not necessary for classification alone. To illustrate this point, consider a particular case with $k=2$ and $Q^*_1$ and $Q^*_2$ are relatively close to each other such that $\overline{Q}=.5\times(Q^*_1+Q^*_2)$ is a good estimation of both $Q^*_1$ and $Q^*_2$. In that case, the model used for the label shift quantification is the singleton containing only $\overline{Q}$, and the estimation of $\beta^*$ is not possible. In that case, estimating $Q^*_1$ or $Q^*_2$ is not as relevant as solving the associated classification problem, wether an observation "comes from" $Q^*_1$ or $Q^*_2$. This observation motivates our study of Setting \ref{case:predictor} in the next section.
\section{Label shift quantification in Setting \ref{case:predictor}}
\label{sec:main}
\subsection{Preliminaries}
Our strategy for the second setting (\ref{case:predictor}) is not really different from the previous one. To improve clarity, we provide a heuristic explanation and establish connection between predictors, label distributions and conditional distributions. Let us use the notation of the introduction briefly where $D_s$ is the source domain, i.e. the distribution of the couple $(X,Y)$ over $\mathscr{X}\times\mathscr{Y}$. The source domain takes the form
\[ D_s(dx,dy)=\alpha^*_y Q^*_y(dx), \]
with $Q^*_1,\dots,Q^*_k\in\mathscr{P}_X$. Let $\mu$ be any $\sigma$-finite (positive) measure dominating $Q^*_1,\dots,Q^*_k$, e.g. $\mu=Q^*_1+\dots+Q^*_k$. We denote by $q^*_i$ the Radon-Nikodym derivative $\mathrm{d}Q^*i/\mathrm{d}\mu$ for all $i\in[k]$. In that case, the Bayes predictor is $f^{\alpha^*}:\mathscr{X}\rightarrow\mathcal{W}_k$, where $f^{\alpha}$ is defined by
\begin{equation}
\label{eq:f_bayes}
f^{\alpha}_i(x) := \frac{\alpha_i q^*_i(x)}{\sum\limits_{j=1}^k \alpha_j q^*_j(x)},
\end{equation}
for all $i\in[k]$ and $x$ in the support of
\begin{equation}
\label{eq:p_alpha}
P_{Q^*,\alpha} := \alpha_1 Q^*_1 + \dots + \alpha_k Q^*_k,
\end{equation}
for all $\alpha\in\mathcal{W}_k$. From (\ref{eq:f_bayes}) and (\ref{eq:p_alpha}), we can deduce
\begin{equation}
\label{eq:f_to_q}
Q_i^*(dx) = (\alpha_i)^{-1} f_i^{\alpha}(x) P_{Q^*,\alpha}(dx),
\end{equation}
for all $i\in[k]$, and all $\alpha$ in
\begin{equation}
\label{eq:W_k_star}
\mathcal{W}_k^*:=\{\beta\in\mathcal{W}_k;\beta_i>0,\forall i\in[k]\}.
\end{equation}
One can check that (\ref{eq:f_to_q}) is valid for any Bayes predictor $f^{\alpha}$ with $\alpha\in\mathcal{W}_k^*$, not only for $f^{\alpha^*}$. This means that knowing $\alpha^*$ is not necessary to have a good predictor. In particular, Section \ref{sec:calibration} indicates that calibration is far more important than the knowledge, or a good estimate, of $\alpha^*$.
\subsection{Our estimator}
\label{sec:rho_est_pred}
Now that we established the link between the conditional distributions and the predictor, we can proceed as in Section \ref{sec:mpe}. Notice that we do not need to know the distributions but only the densities to construct $\rho$-estimators in Section \ref{sec:rho_est}. Therefore, in the ideal scenario where we have access to the Bayes estimator $f^{\alpha}$, with $\alpha$ in $\mathcal{W}_k^*$ defined by (\ref{eq:W_k_star}), we can construct our estimator using densities $q^*_i= (\alpha_i)^{-1} f^{\alpha}_i$ for all $i\in[k]$, even if we do not know the conditional distributions $Q^*_1,\dots,Q^*_k$. We simply extend this approach by plugging in the predictor $f$ and the label probability $\alpha\in\mathcal{W}_k$.\\
Given a predictor $f:\mathscr{X}\rightarrow \mathcal{W}_k$ and weights $(\alpha_i)_i$ with $\alpha_i>0$ for all $i\in[k]$, we consider a mixture model with fixed emission distributions/densities. We define the countable class of functions $\mathcal{M}(f,\alpha)$ by
\[ \mathcal{M}(f,\alpha) := \left\{ x\in\mathscr{X} \mapsto \sum\limits_{i=1}^k \beta_i \alpha_i^{-1} f_i(x); \beta \in\mathcal{W}_k\cap\mathbb{Q}^k \right\}. \]
Although it is not necessary for the construction of our estimator, we still need to associate a class of probability distributions to $\mathcal{M}(f,\alpha)$ to quantify its performance. We define the class of measures
\[ \mathscr{P}(f,\alpha) := \left\{ \text{positive }\sigma\text{-finite measure  } \mu \text{  on  } (\mathscr{X},\mathcal{X}) \text{  such that  } \mu(f_i) = \alpha_i, \forall i\in[k] \right\}, \]
where $\mu(f_i)$ denotes the integral of $f_i$ with respect to $\mu$. For $\mu$ in $\mathscr{P}(f,\alpha)$, we define the model
\begin{equation}
\label{eq:overline_M}
\overline{\mathscr{M}}(f,\alpha,\mu) := \left\{ \sum\limits_{i=1}^k \beta_i \alpha_i^{-1} f_i \cdot \mathrm{d}\mu; \beta\in\mathcal{W}_k \right\}.
\end{equation}
As discussed in Section \ref{sec:mpe}, the linear independence of conditional distributions is necessary in order to correctly define our estimator. Therefore, we define the class of measures
\[ \mathscr{P}^*(f,\alpha) := \left\{ \mu \in \mathscr{P}(f,\alpha); \begin{array}{c}
\text{  distributions  } f_1\cdot \mathrm{d}\mu,\dots, f_k\cdot\mathrm{d}\mu\\
\text{  are linearly independent  }
\end{array}
\right\}, \]
and we make the following assumption.
\begin{assumption}
\label{hyp:nu}
The class of distributions $\mathscr{P^*}(f,\alpha)$ is not empty. 
\end{assumption}
Notice that it is the first assumption we make on the predictor $f$ (and $\alpha$), in particular, we did not assume that it is calibrated. Calibration and its role in label shift quantification is discussed in Section \ref{sec:calibration}.\\
Our estimator is defined as follows. Let $\mu$ be in $\mathscr{P}(f,\alpha)$ and $\hat{P}=\hat{P}\left(\mathbf{X},\mathcal{M}(f,\alpha)\right)\in\overline{\mathscr{M}}(f,\alpha,\mu)$ be a $\rho$-estimator, as defined by (\ref{eq:rho_est}). We denote by $\hat{\beta}$ any element of $\mathcal{W}_k$ such that
\begin{equation}
\label{eq:hat_beta_predictor}
\hat{P}=\sum\limits_{i=1}^k \hat{\beta}_i \alpha_i^{-1} f_i \cdot \mathrm{d}\mu.
\end{equation}
Note that if $\mu$ belongs to $\mathscr{P}^*(f,\alpha)$, this element is unique and it does not depend on $\mu$. Then we say that $\hat{\beta}$ is a $\rho$-estimator. Note that assuming $Q^*_1,\dots,Q^*_k$ to be linearly independent implies that $\mathscr{P}^*(f^{\alpha},\alpha)$ is non-empty as it contains $P_{Q^*,\alpha}$. In that case, we have $\overline{\mathscr{M}}(f,\alpha,\mu)=\overline{\mathscr{M}}_{mix}(Q_1,\dots,Q_k)$ with $Q_i=\alpha_i^{-1} f_i \cdot \mathrm{d}\mu$ for all $i\in[k]$. This means we have similar results to those in Section \ref{sec:mpe}, and, in particular, the connection with the MLE.
\begin{prop}
\label{prop:mle_predictor}
When it exists, the Maximum Likelihood Estimator $\hat{\beta}_{MLE}$ given by
\begin{equation}
\label{eq:mle_b}
\hat{\beta}_{MLE} \in \argmax_{\beta\in\mathcal{W}_k} \sum\limits_{i=1}^n \log\left( \sum\limits_{j=1}^k \beta_j \alpha_j^{-1} f_j(X_i) \right)
\end{equation}
is a $\rho$-estimator with respect to $\overline{\mathscr{M}}(f,\alpha,\mu)$ for any $\mu\in\mathscr{P}^*(f,\alpha)$.
\end{prop}
This result is proven in Section \ref{sec:proof_prop_mle_predictor}. This result implies that all the results we give in this section are also valid for maximum likelihood label shift (MLLS) -- see \cite{garg} or \cite{alexandari}. 
\subsection{Results}
We assume that Assumption \ref{hyp:nu} holds in all this section. From now on, $\hat{\beta}$ will denote either the estimator given by (\ref{eq:hat_beta_predictor}) or the MLE defined by (\ref{eq:mle_b}), where $\mu$ is any element of $\mathscr{P}^*(f,\alpha)$.
\begin{theorem}
\label{th:predictor}
\begin{itemize}
\item  There is a positive constant $C(f,\alpha,\mu)$ depending only on $f,\alpha$ and $\mu$, such that for all $\overline{\beta}\in\mathcal{W}_k$ and for all $\xi>0$,
\[ C(f,\alpha,\mu) ||\hat{\beta}-\overline{\beta}||_1^2 \leq n^{-1} \sum\limits_{j=1}^n h^2\left( P_j, \sum\limits_{i=1}^k \overline{\beta}_i \alpha_i^{-1} f_i \cdot \mathrm{d}\mu \right) + \frac{k \log (n/k) + \xi}{n}, \]
with probability at least $1-e^{-\xi}$.
\item Let $\overline{Q}_1,\dots,\overline{Q}_k$ be linearly independent distributions in $\mathscr{P}_X$. There is a positive constant $C(\overline{\mathbf{Q}})$ depending only on $\overline{Q}_1,\dots,\overline{Q}_k$ such that for all $\overline{\beta}\in\mathcal{W}_k$ and all $\xi>0$,
\begin{align*}
C(\overline{\mathbf{Q}}) ||\hat{\beta}-\overline{\beta}||_1^2 &\leq n^{-1} \sum\limits_{j=1}^n h^2\left( P_j, \sum\limits_{i=1}^k \overline{\beta}_i \overline{Q}_i \right) + \max_{i\in[k]} h^2\left( \alpha_i^{-1}f_i\cdot \mathrm{d}\mu, \overline{Q}_i \right)\\
&+ \frac{k \log (n/k) + \xi}{n},
\end{align*}
with probability at least $1-e^{-\xi}$. 
\end{itemize}
\end{theorem}
The proof can be found in Section \ref{sec:proof_th_predictor}. As Theorem \ref{th:conditional}, this result is general but is not very interesting unless we make assumptions on the distributions $P_1,\dots,P_n$. The quantity
\[ \sum\limits_{j=1}^n h^2\left( P_j, \sum\limits_{i=1}^k \overline{\beta}_i \alpha_i^{-1} f_i \cdot \mathrm{d}\mu \right) \]
quantifies the distance from our model to the true distributions of the observations. We first put ourselves in the context of Assumption \ref{hyp:iid_mix} to understand the influence of $f$ and $\alpha$ on this quantity. We introduce the notion of confusion matrix.
\begin{definition}{Confusion matrix}\\
Under Assumption \ref{hyp:iid_mix}, for a predictor $f:\mathscr{X}\rightarrow \mathcal{W}_k$ we denote by $M(f)$ the \emph{confusion matrix} in $\mathbb{R}^{k\times k}$ defined by
\[ M(f)_{ij} := Q^*_j(f_i) = \int f_i(x) Q^*_j(dx), \]
for all $i,j\in[k]$.    
\end{definition}
One can see that the confusion matrix is a stochastic matrix since $f$ takes values in $\mathcal{W}_k$ and $Q^*_1,\dots,Q^*_k$ are probability distributions, i.e. $M(f)_{1j}+\dots+M(f)_{kj}=1$ for all $j\in[k]$.
\begin{assumption}
\label{hyp:q_sigma}
\begin{enumerate}
\item The measures $f_1 \cdot \mathrm{d}Q^*_{\Sigma},\dots,f_k\cdot \mathrm{d}Q^*_{\Sigma}$ are linearly independent, where $Q^*_{\Sigma}=Q^*_1+\dots+Q^*_k$.
\item There is $\gamma \in \mathbb{R}^k$ such that $\sum\limits_{i=1} \gamma_i Q^*_i \in\mathscr{P}^*(f,\alpha)$ with $\gamma_i>0$ for all $i\in[k]$.
\end{enumerate}
\end{assumption}
This assumption might appear unusual but the next result indicates that it is rather standard and it is satisfied for $f=f^{\alpha}$ and $\alpha=\alpha$ in particular. The second claim is weaker than the standard calibration assumptions -- see Section \ref{sec:calibration} -- given the first claim. Lemma \ref{lem:marginal_calibration} below shows that calibration implies $\alpha = M(f)\alpha$ and one can easily deduce that $\sum_{i=1}^k \alpha_i Q^*_i$ is a probability distribution.
\begin{prop}
\label{prop:confusion_invertible}
\begin{itemize}
    \item If $f_1\cdot \mathrm{d}Q^*_{\Sigma},\dots,f_k\cdot \mathrm{d}Q^*_{\Sigma}$ are linearly independent, then the confusion matrix $M(f)$ is invertible. If $\gamma:=M(f)^{-1}\alpha$ has positive coordinates, then $\sum\limits_{i=1}^k \gamma_i Q^*_i \in \mathscr{P}^*(f,\alpha)$.
    \item Consider the ideal case $f=f^{\alpha}$. We have $\alpha = M(f^{\alpha})\alpha$, and therefore $\sum\limits_{i=1}^k \alpha_i Q^*_i \in \mathscr{P}^*(f^{\alpha},\alpha)$. We have
    \[ Q^*_1,\dots,Q^*_k \text{ are linearly independent} \Leftrightarrow f^{\alpha}_1\cdot \mathrm{d}Q^*_{\Sigma},\dots,f^{\alpha}_k\cdot \mathrm{d}Q^*_{\Sigma} \text{ are linearly independent}. \]
\end{itemize}
\end{prop}
This result is proven in Section \ref{sec:proof_prop_confusion_invertible}. It shows that Assumption \ref{hyp:q_sigma} is not too restrictive and is satisfied for $f=f^{\alpha}$ and $\alpha=\alpha$ in particular. In general, checking that $\gamma=M(f)^{-1}\alpha$ has positive coordinates is not straightforward but it is not a very strong assumption. Since $\alpha$ is a fixed point of $M(f^{\alpha})$, it is natural to expect that $M(f)^{-1}\alpha$ has positive coordinates when $f$ is close to $f^{\alpha}$, which should be the case for a good predictor. Moreover, the divergence of $\gamma$ from $\alpha$ is related to the divergence of $f$ from $f^{\alpha}$ through the confusion matrices,
\[ ||\alpha - \gamma||_1 = ||\alpha - M(f)^{-1}\alpha||_1 = ||(I_k - M(f)^{-1} M(f^{\alpha}))\alpha||_1 \leq ||I_k-M(f)^{-1} M(f^{\alpha})||_{1,1}, \]
where $||\cdot||_{1,1}$ is the operator norm associated to the $\ell_1$-norm.\par
We have the following deviation inequality.
\begin{coro}
\label{coro:predictor}
Under Assumptions \ref{hyp:iid_mix} and \ref{hyp:q_sigma}, for all $\xi>0$ we have
\[ C(f,\alpha,\mu) ||\beta^*-\hat{\beta}||_1^2 \leq \max\limits_{i\in[k]} (\beta^*_i/\alpha_i) \left( \mathbb{E}_{P_{Q^*,\alpha}}\left[ ||f^{\alpha} - f||_1 \right] + || \alpha - \gamma ||_1 \right) + \frac{k \log (n/k) + \xi}{n}, \]
with probability at least $1-e^{-\xi}$, where $f^{\alpha}$ is given by (\ref{eq:f_bayes}), $P_{Q^*,\alpha}$ is given by (\ref{eq:p_alpha}) and $\gamma$ is given in Assumption \ref{hyp:q_sigma}.
\end{coro}
This result is proven in Section \ref{sec:proof_coro_predictor}. As in Setting \ref{case:1}, in the well-specified case, i.e. for $f=f^{\alpha}$, we have
\begin{equation}
\label{eq:well_specified_predictor}
C(\mathbf{Q}^*) ||\beta^*-\hat{\beta}||_1^2 \leq \frac{k \log (n/k) + \xi}{n},
\end{equation}
with probability at least $1-e^{-\xi}$, for all $\xi>0$. In that case, the constant depends on $Q^*_1,\dots,Q^*_k$ or equivalently on $f^{\alpha},\alpha$ and $P_{Q^*,\alpha}$. We have a bound on the convergence rate of our estimator of order $(n/k)^{-1/2}\log^{1/2}(n/k)$.\\
This result is similar to \cite[Theorem 3]{garg} with a few differences. They give a bound on the $\ell_2$-distance $||w-\hat{w}||_2$ where $w$ is the shift ratio defined by $w_i=\beta^*_i/\alpha_i$ for all $i\in[k]$, and $\hat{w}_i=\hat{\beta}_i/\alpha_i$ for all $i\in[k]$. It is not clear if this is the reason but there are also differences in the constants appearing  in both results. Their result includes a term $\sigma_{f}^{-1}$, the minimal eigenvalue of $k\times k$ matrix, which seems to combine both the roles of $C(\mathbf{Q}^*)$ and $k$ in (\ref{eq:well_specified_predictor}). One major differences with \cite[Theorem 3]{garg} is the fact that they weaken the assumption $f=f^{\alpha}$ to a calibration assumption. This is addressed in Section \ref{sec:calibration} with Theorem \ref{th:calibration} with weaker assumptions since it does include anything like Condition 1 \cite{garg}. However, we illustrate below the robustness properties of the estimator which is not the case in \cite{garg}.\par
As in Setting \ref{case:1}, the estimator $\hat{\beta}$ possesses robustness properties. Corollary \ref{coro:predictor} already included some robustness to misspecification, indicating that the performance of our estimator is not significantly worse as long as the quantity
\[ \max_{i\in[k]} (\beta^*_i/\alpha_i)\left( \mathbb{E}_{P_{Q^*,\alpha}}\left[ ||f^{\alpha}-f||_1 \right] + ||\alpha-\gamma||_1 \right) \]
is of order not greater than $n^{-1} k \log(n/k)$. Following the remark after Proposition \ref{prop:confusion_invertible}, $||\alpha - \gamma||_1$ can be seen as a measure of the divergence of $f$ from the Bayes predictor $f^{\alpha}$ -- through the confusion matrix in a sense.\\
We can also consider more specific cases of misspecification, such as contamination or the presence of outliers. To simplify the analysis, we assume that $f=f^{\alpha}$ which means that we satisfy Assumption \ref{hyp:q_sigma} and the first inequality of Theorem \ref{th:predictor} becomes
\[ C(\mathbf{Q}^*) ||\hat{\beta} - \overline{\beta}||_1^2 \leq \frac{1}{n} \sum\limits_{j=1}^n \left( P_j, \sum\limits_{i=1}^k \overline{\beta}_i Q^*_i \right) + \frac{k\log(n/k) + \xi}{n}, \]
which holds with probability at least $1-e^{-\xi}$ for all $\xi>0$. If we consider the case of contamination, as in (\ref{eq:p_contamination}), for all $\xi>0$, we have
\[ C(\mathbf{Q}^*) ||\hat{\beta}-\beta^*||_1^2 \leq \lambda_0 + \frac{k \log (n/k) + \xi}{n}, \]
with probability at least $1-e^{\xi}$, where $\beta^*_i=\lambda_i/(1-\lambda_0)$ for all $i\in[k]$. As long as $\lambda_0$ is small compared to $n^{-1} k \log (n/k)$ the performance of our estimator is not significantly worse than in the ideal setting. Similarly, if we consider the potential presence of outliers as in (\ref{eq:p_outliers}), for all $\xi>0$, we have
\[ C(\mathbf{Q}^*) ||\beta^*-\hat{\beta}||_1^2 \leq \frac{|I|}{n} + \frac{k \log (n/k)   + \xi}{n},
\]
with probability at least $1-e^{-\xi}$. As long as the proportion of outliers $|I|/n$ is small compared to $n^{-1} k\log (n/k)$ the performance of the estimator is not significantly worse.\par
Although we considered separately misspecification, outliers and contamination, it is possible to combine those different cases and obtain a bound for the estimator. Generally speaking, ss long as the departure from the ideal setting is not too important, it does not significantly affect the performance of our estimator.\\
We can see robustness. It backs up the numerical study of Saerens \emph{et al}. (Section 4), from which they draw the conclusion that 'the EM algorithm appeared to be more robust than the
confusion matrix method'. The interpretation of Theorem \ref{th:predictor} and Corollary \ref{coro:predictor} seems to confirm the intuition of \cite{alexandari}.
\subsection{Calibration}
\label{sec:calibration}
Until now, we briefly mentioned calibration but did not really consider it. Calibration is a crucial point of MLLS (see \cite{alexandari,garg,kumar,vaicenavicius}). Several papers, such as \cite{alexandari} and \cite{garg}, showed empirically and theoretically that MLLS will outperform BBSE, but only if the predictor $f$ is well calibrated. We would like to reach similar conclusions but in order to utilize the concept of calibration, we need a definition that is suitable for our context. In particular, the definition one can find in \cite{garg} is not suitable for our very formal context with general spaces $\mathscr{X}$. Instead, we use the definitions given in \cite{vaicenavicius} which are more adapted here.\par
\subsubsection{Definitions and preliminary results}
We will use the following notation. For a finite set $I$, we denote by $\mathcal{W}_I$ the simplex of dimension $|I|$ defined by
\[ \mathcal{W}_I := \left\{ (w_i)_{i\in I} \in [0,1]^{|I|} : \sum\limits_{i\in I} w_i = 1 \right\}. \]
As we did earlier, we will write $\mathcal{W}_k=\mathcal{W}_{[k]}$ for all $k\geq 1$.
Let $(\mathscr{X},\mathcal{X})$ be a measurable space. Let $\Pi$ be a distribution on $(\mathscr{X},\mathcal{X})\otimes (I,\mathcal{P}(I))$, where $\mathcal{P}(I)$ is the set of all subsets of $I$. We will just say that $\Pi$ is a distribution on $I\times\mathscr{X}$ for simplicity, from now on.
\begin{definition}{(Calibration \cite{vaicenavicius})}\\
\begin{itemize}
    \item We say that a predictor $f$ is \emph{marginally} calibrated, with respect to a distributon $\Pi$ on $I\times\mathscr{X}$, if for all $i\in I$,
\[ f_i(X) = \mathbb{E}\left[ \mathbbm{1}_{Y=i}| f_i(X) \right], \]
where $(X,Y)\sim\Pi$.
    \item We say that $f:\mathscr{X}\rightarrow\mathcal{W}_I$ is \textbf{canonically calibrated}, with respect to a distribution $\Pi$ on $I\times \mathscr{X}$, if for all $i\in I$,
\[ f_i(X) = \mathbb{E}_{\Pi}[\mathbbm{1}_{Y=i} | f(X) ], \]
where $(X,Y)\sim\Pi$.
\end{itemize}
\end{definition}
One can easily see that any canonically calibrated predictor is marginally calibrated. We mainly focus on canonical calibration from now on. The following result shows how canonical calibration is preserved under product operations.
\begin{lemma}
\label{lem:calibration_product}
Let $I_1$ and $I_2$ be two finite sets. Let $(\mathscr{X}_1,\mathcal{X}_1)$ and $(\mathscr{X}_2,\mathcal{X}_2)$ be two measurable spaces. If $f_1:\mathscr{X}_1\rightarrow\mathcal{W}_{I_1}$ is canonically calibrated with respect to $\Pi_1$ and $f_2:\mathscr{X}_2\rightarrow\mathcal{W}_{I_2}$ is canonically calibrated with respect to $\Pi_2$, then
\[ g: 
\begin{cases}
\mathscr{X}_1\times\mathscr{X}_2\rightarrow\mathcal{W}_{I_1\times I_2}\\
(x_1,x_2)\mapsto (f_{i_1}(x_1)f_{i_2}(x_2))_{i_1,i_2} 
\end{cases}
\]
is canonically calibrated with respect to $\Pi$, given by
\[ (X_1,Y_1,X_2,Y_2) \sim \Pi_1\otimes \Pi_2 \Leftrightarrow ((X_1,X_2),(Y_1,Y_2))\sim \Pi. \] 
\end{lemma}
This resulted is proven in Section \ref{sec:proof_lem_calibration_product}. By induction, we can extend this to any $n\geq 1$. In particular, we have the following result.
\begin{coro}
\label{coro:calibration_product}
Let $I$ be a finite set and $(\mathscr{X},\mathcal{X})$ be a measurable space. Let $n$ be any integer larger than 1. If $f:\mathscr{X}\rightarrow\mathcal{W}_{I}$ is canonically calibrated with respect to $\Pi$, then
\[ g :
\begin{cases}
\mathscr{X}^n \rightarrow \mathcal{W}_{I^n}\\
(x_1,\dots,x_n)\mapsto (f_{i_1}(x_1) \dots f_{i_n}(x_n))_{i_1,\dots,i_n} 
\end{cases}
\]
is canonically calibrated with respect to $\Pi^{(n)}$ given by
\[ (X_1,Y_1,\dots,X_n,Y_n) \sim \Pi^{\otimes n}\Leftrightarrow ((X_1,\dots,X_n),(Y_1,\dots,Y_n))\sim \Pi^{(n)}. \] 
\end{coro}
The next lemma is probably one of the most important of this article. It gives an interpretation of canonical calibration that shows it is as if we had access to the Bayes predictor, up to the linear independence assumption. Let $\Pi$ be a distribution on $I\times\mathscr{X}$. We define
\[ (w_i)_{i\in I}:= \left( \Pi(\mathscr{X}\times\{i\}) \right)_{i\in I} \in\mathcal{W}_I, \]
the distribution of $Y$ under $\Pi$, and
\[ R_i(dx) := \Pi(\{i\},dx), \]
the distribution of $X$ conditionally to $Y=i$ under $\Pi$, for all $i\in I$.
\begin{lemma}
\label{lem:canonical_calibration_exp}
Let $f:\mathscr{X}\rightarrow \mathcal{W}_I$ be canonically calibrated with respect to $\Pi$. Assume $w_i>0$ for all $i\in I$. For $v\in\mathcal{W}_I$, we define the following distributions over ($\mathscr{X},\mathcal{X}$)
\[ P_{R,v}(dx) := \sum\limits_{i\in I} v_i R_i(dx) \quad \text{  and  } \quad P_{f,v}(dx) := \sum\limits_{i\in I} v_i w_i^{-1} f_i(x) P_{R,w}(dx). \]
For all measurable functions $\phi:\mathcal{W}_I\rightarrow\mathbb{R}$ and all $v \in\mathcal{W}_I$ we have
\[ \mathbb{E}_{P_{R,v}}\left[ \phi(f(X)) \right] = \mathbb{E}_{P_{f,v}}\left[ \phi(f(X))\right]. \]
\end{lemma}
This result is proven in Section \ref{sec:proof_lem_canonical_calibration_exp}. The distribution $P_{R,v}$ is the same distribution as $\Pi$, up to label shift whereas $P_{f,v}$ is the distribution with a second order shift, as we have a shift on the distribution of $Y$ and a shift on the distribution of $X$ when we use the predictor $f$ to construct a distribution as if we had access to the Bayes predictor -- i.e. $P_{f^{w},v}=P_{R,v}$ if $f^{w}$ is the Bayes predictor.
\subsubsection{Consequences of marginal calibration}
To understand better the role of calibration, we forget about robustness considerations and work under Assumption \ref{hyp:iid_mix}. For $\alpha$ in $\mathcal{W}_k$, we denote by $\Pi_{\alpha}$ the probability distribution over $\mathscr{X}\times [k]$ defined by
\begin{equation}
\label{eq:pi_alpha}
\Pi_{\alpha}(A\times\{i\})=\alpha_i Q^*_i(A),
\end{equation}
for all $i\in[k]$ and all measurable sets $A\subset \mathcal{X}$. One can check that the Bayes predictor $f^{\alpha}$ given by (\ref{eq:f_bayes}) is (canonically) calibrated with respect to $\Pi_{\alpha}$, as well as the constant predictor $f:x\mapsto \alpha$. The next result connects the notion of calibration with the matrix confusion.
\begin{lemma}
\label{lem:marginal_calibration}
Under Assumption \ref{hyp:iid_mix}, if $f$ is marginally calibrated with respect to $\Pi_{\alpha}$, we have $\alpha = M(f) \alpha$ and therefore $P_{Q^*,\alpha}\in\mathscr{P}(f,\alpha)$, where $P_{Q^*,\alpha}$ is given by (\ref{eq:p_alpha}).
\end{lemma}
Particularly, the second claim of Assumption \ref{hyp:q_sigma} is satisfied if $f$ is marginally calibrated with respect to $\Pi_{\alpha}$. We have the following result .
\begin{coro}
\label{coro:predictor_calibration}
Under Assumptions \ref{hyp:iid_mix} and \ref{hyp:q_sigma}, if $f$ is marginally calibrated with respect to $\Pi_{\alpha}$, for all $\xi>0$,
\begin{equation}
\label{eq:coro_predictor_calibration}
C(\mathbf{Q}^*) ||\beta^*-\hat{\beta}||_1^2 \leq \mathbb{E}_{P_{Q^*,\alpha}}\left[ \left| \sum\limits_{i=1}^k \beta^*_i \alpha_i^{-1} (f_i^{\alpha} - f_i)(X) \right| \right] + \frac{k \log(n/k) + \xi}{n},
\end{equation}
with probability at least $1-e^{-\xi}$, where $C(\mathbf{Q}^*)$ is a positive constant only depending on $Q^*_1,\dots,Q^*_k$.
\end{coro}
This result is proven in Section \ref{sec:proof_coro_predictor_calibration}. It is an improvement of Corollary \ref{coro:predictor} as we have
\[ \mathbb{E}_{P_{Q^*,\alpha}}\left[ \left| \sum\limits_{i=1}^k \beta^*_i \alpha_i^{-1} (f_i^{\alpha} - f_i)(X) \right| \right] \leq \max_{i\in[k]} (\beta^*_i/\alpha_i)  \mathbb{E}_{P_{Q^*,\alpha}}\left[ ||f^{\alpha}-f||_1 \right]. \]
Marginal calibration allows us to remove the term $||\alpha-\gamma||$ from the bound but it does not seem to be enough to obtain better results such as consistency for instance.
\subsubsection{Consequences of canonical calibration}
The next results show that canonical calibration is a sufficient condition to obtain a consistent estimator in the context of label shift quantification. Assume that $f$ is canonically calibrated with respect to $\Pi_{\alpha}$. In that case, $P_{Q^*,\alpha}$ belongs to $\mathscr{P}(f,\alpha)$ and the distribution set given by (\ref{eq:overline_M}) can be written as $\overline{\mathscr{M}}(f,\alpha,P_{Q^*,\alpha})= \left\{ P_{f,\alpha,\beta}; \beta\in\mathcal{W}_k \right\}$,
where
\begin{equation}
\label{eq:p_f_beta}
P_{f,\alpha,\beta} := \sum\limits_{i=1}^k \beta_i \alpha_i^{-1} f_i(x) P_{Q^*,\alpha}(dx),
\end{equation}
for all $\beta\in\mathcal{W}_k$.
\begin{prop}
\label{prop:kl_min_canon}
If $f$ is canonically calibrated with respect to $\Pi_{\alpha}$,
\[ \argmax_{\beta\in\mathcal{W}_k} \mathbb{E}_{P^*}\left[ \log\left( \sum\limits_{i=1}^k \beta_i \alpha_i^{-1} f_i(X) \right) \right] = \argmin_{\beta\in\mathcal{W}_k} \mathbf{K}( P_{f,\beta^*}||P_{f,\beta}), \]
where $\mathbf{K}$ denotes the Kullback-Leibler divergence. 
\end{prop}
This result is proven in Section \ref{sec:proof_prop_kl_min_canon}. It shows that if $f$ is canonically calibrated with respect to $\Pi_{\alpha}$, the KL divergence minimizer is the same as the likelihood maximizer. In particular, if $P_{Q^*,\alpha}$ belongs to $\mathscr{P}^*(f,\alpha)$, the KL divergence minimizer is $\beta^*$ and therefore $\beta^*$ is the likelihood maximizer. This result is not new as it is a formal reformulation of Lemma 2 in \cite{garg}. The interest behind this result is the tools developed to prove it, which gives a deeper understanding of the role of calibration in label shift quantification -- particularly Lemma \ref{lem:calibration_product}.\par
Finally, we are able to prove that canonical calibration and the linear independence assumption are sufficient to obtain a consistent estimator with the desired convergence rate.
\begin{theorem}
\label{th:calibration}
Let Assumption \ref{hyp:iid_mix} hold. Let $f$ be canonically calibrated with respect to $\Pi_{\alpha}$ given by (\ref{eq:pi_alpha}). For all $\xi>0$, with probability at least $1-e^{-\xi}$, we have
\begin{equation}
\label{eq:th_calibration_1}
h^2\left( P_{f,\alpha,\beta}, P_{f,\alpha,\hat{\beta}} \right) \leq c_2 \frac{k\log(n/k)}{n} + c_3 \frac{\xi}{n}.
\end{equation}
In particular, if $P_{Q^*,\alpha}\in\mathscr{P}^*(f,\alpha)$ -- i.e. $f_1\cdot \mathrm{d}P_{Q^*,\alpha},\dots,f_k\cdot \mathrm{d}P_{Q^*,\alpha}$ are linearly independent -- for all $\xi>0$,
\begin{equation}
\label{eq:th_calibration_2}
C(\mathbf{Q}^*,f) ||\beta^*-\hat{\beta}||^2_1 \leq \frac{k\log(n/k) + \xi}{n},
\end{equation}
with probability at least $1-e^{-\xi}$, where $C(\mathbf{Q}^*,f)$ is a positive constant only depending on $Q^*_1,\dots,Q^*_k$ and $f$.
\end{theorem}
This result is proven in Section \ref{sec:proof_th_calibration}. It shows that if $f$ is canonically calibrated with respect to $\Pi_{\alpha}$, we have the same convergence rate as in the ideal setting where we dispose of the Bayes predictor. However, we still pay a price for using a different predictor with a constant $C(\mathbf{Q}^*,f)$ that can be much smaller than $C(\mathbf{Q}^*)$ and therefore a worse convergence rate.
\section{Conclusion}
This paper provides theoretical guarantees for label shift quantification using off-the-shelf conditional distributions or predictors. Specifically, we establish convergence rate bounds in the well-specified case and demonstrate robustness to outliers and contamination for the proposed method, which includes Maximum Likelihood Label Shift. Our findings support and extend the numerical study of Saerens \emph{et al}. (Section 4), confirming the robustness properties of MLLS and further strengthening the theoretical foundation for their use. This work complements the contributions of \cite{saerens}, \cite{alexandari}, and \cite{garg}, offering a comprehensive perspective of MLLS in label shift estimation. Finally, we introduce a formalism for calibration, which gives a better understanding of its implications -- see Lemma \ref{lem:canonical_calibration_exp} for instance.
\section*{\fontsize{12}{15}\selectfont Acknowledgements}
The author has been partially supported by ANR-21-CE23-0035 (ASCAI).
\bibliography{bibli}
%
\begin{appendix}
\section{Proofs of Section \ref{sec:mpe}}
The following lemma is central and allows us to prove Theorems \ref{th:conditional} and \ref{th:predictor}.
\begin{lemma}
\label{lem:aux_2}
Let $n\geq e \times k$. Let $\hat{P}$ be the $\rho$-estimator on $\mathbf{X}=(X_1,\dots,X_n)$ and $\mathcal{M}_{mix}(q_1,\dots,q_k)$. If $X_1,\dots,X_n$ are independent with distribution $P_1,\dots,P_n$, for all $\xi>0$ and all $\overline{w}\in\mathcal{W}_k$, we have
\[ \sum\limits_{j=1}^n h^2\left( P_j, \hat{P} \right) \leq c_1 \sum\limits_{j=1}^n h^2\left( P_j, \sum\limits_{i=1}^k \overline{w}_i Q_i \right) +  c_2 k \log(n/k) + c_3 \xi, \]
with probability at least $1-e^{-\xi}$, where $c_1=150$, $c_2=2.1 \times 10^6$ and $c_3=5014$. 
\end{lemma}
This result is proven in Section \ref{sec:proof_lem_aux2}.
\subsection{Proof of Proposition \ref{prop:mle}}
\label{sec:proof_prop_mle}
It is a direct consequence of Corollary 1 of \cite{baraudrevisited} since
\[ \overline{\mathcal{M}}_{mix}(q_1,\dots,q_k) := \left\{ \sum\limits_{i=1}^k \beta_i q_i; \beta \in \mathcal{W}_k  \right\} \]
is a convex set of densities.
\subsection{ Proof of Lemma \ref{lem:aux_1}}
\label{sec:proof_lem_aux1}
Let $F_1,\dots,F_k$ be distributions in $\mathscr{P}_X$. Let $\lambda$ be a $\sigma$-finite measure dominating $F_1,\dots,F_K$. One can always take $\lambda=F_1+\dots+F_K$. Let $f_1,\dots,f_k$ be the respective density functions of $F_1,\dots,F_k$ with respect to $\lambda$. For $\beta\in\mathcal{W}_K$, we write
\[ P_{\beta} = \beta_1 F_1 + \dots + \beta_K F_K  \text{  and  } p_{\beta} = w_1 f_1 + \dots + w_K f_K. \]
Fix two elements $\beta\neq \beta'\in\mathcal{W}_K$. Let us define
\[ \delta = \sum\limits_{i\in I} \beta_i-\beta'_i = \sum\limits_{i\in I^c} \beta'_i-\beta_i = d_{TV}(\beta,\beta') = \frac{1}{2} ||\beta-\beta'||_1, \]
where $I=\{ i\in[k] : \beta_i \geq \beta'_i \}$ and $I^c=[k]\backslash I$. We also define
\[ Q_I = \sum\limits_{i\in I} a_i F_i,  \text{  and  } Q_{I^c} = \sum\limits_{i\in I^c} b_i F_i, \]
where $a = \left(\frac{\beta_i-\beta'_i}{\delta}\right)_{i\in I}\in\mathcal{W}_{|I|}$ and $b = \left(\frac{\beta'_i-\beta_i}{\delta}\right)_{i\in [k]\backslash I}\in\mathcal{W}_{|I^c|}$, and the associated densities $q_I = \sum\limits_{i\in I} a_i f_i$ and $q_{I^c} = \sum\limits_{i\in I^c} b_i f_i$. One can easily check that
\[ P_\beta - P_{\beta'} = \delta \left( Q_I - Q_{I^c} \right). \]
Therefore, we have
\[ d_{TV}\left(P_{\beta},P_{\beta'}\right) = \delta d_{TV}(Q_I,Q_{I^c}). \]
We can conclude with the classical inequality
\[ \sqrt{2} h(P,Q)\geq d_{TV}(P,Q), \]
for all $P,Q\in\mathscr{P}_X$.
\subsection{Proof of Theorem \ref{th:conditional}}
\label{sec:proof_th_conditional}
Let $\overline{\beta}$ be in $\mathcal{W}_k$. As a direct consequence of Lemma \ref{lem:aux_2}, there is an event $\Omega_{\xi}$ of probability $1-e^{-\xi}$ such that on $\Omega_{\xi}$, we have
\begin{align*}
\sum\limits_{j=1}^n h^2\left( P_j, \sum\limits_{i=1}^k \hat{\beta}_i Q_i \right) &\leq c_1 \sum\limits_{j=1}^n h^2\left( P_j, \sum\limits_{i=1}^k \overline{\beta} Q_i \right) + c_2 k\log (n/k) + c_3 \xi.
\end{align*}
Using
\begin{align*}
h^2\left( \sum\limits_{i=1}^k \overline{\beta}_i Q_i, \sum\limits_{i=1}^k \hat{\beta}_i Q_i \right) \leq \frac{2}{n} \sum\limits_{j=1}^n h^2\left( \sum\limits_{i=1}^k \overline{\beta}_i Q_i, P_j \right) + \frac{2}{n} \sum\limits_{j=1}^n h^2\left( P_j, \sum\limits_{i=1}^k \hat{\beta}_i Q_i \right), 
\end{align*}
we get
\[ h^2\left( \sum\limits_{i=1}^k \overline{\beta}_i Q_i, \sum\limits_{i=1}^k \hat{\beta}_i Q_i \right) \leq \frac{2(1+c_1)}{n} \sum\limits_{j=1}^n h^2\left( P_j, \sum\limits_{i=1}^k \overline{\beta} Q_i \right) + \frac{2 c_2 k\log (n/k)}{n} + \frac{2c_3\xi}{n}, \]
on $\Omega_{\xi}$ and Lemma \ref{lem:aux_1} allows us to conclude. We can also have the constant depending on any distributions $\overline{Q}_1,\dots,\overline{Q}_k$ that are linearly independent. On $\Omega_{\xi}$, we have
\begin{align*}
h^2\left( \sum\limits_{i=1}^k \overline{\beta}_i \overline{Q}_i, \sum\limits_{i=1}^k \hat{\beta}_i \overline{Q}_i \right) &\leq \frac{3}{n} \sum\limits_{j=1}^n h^2\left( \sum\limits_{i=1}^k \overline{\beta}_i \overline{Q}_i, P_j \right) + \frac{3}{n} \sum\limits_{j=1}^n h^2\left( P_j, \sum\limits_{i=1}^k \hat{\beta}_i Q_i \right)\\
&+ 3 h^2\left( \sum\limits_{i=1}^k \hat{\beta}_i Q_i, \sum\limits_{i=1}^k \hat{\beta}_i \overline{Q}_i \right)\\
&\leq 3 \max_{1\leq i\leq k} h^2\left( Q_i, \overline{Q}_i \right) + \frac{3}{n} \sum\limits_{j=1}^n h^2\left( \sum\limits_{i=1}^k \overline{\beta}_i \overline{Q}_i, P_j \right)\\
&+ \frac{3 c_1}{n} \sum\limits_{j=1}^n h^2\left( P_j, \sum\limits_{i=1}^k \overline{\beta}_i Q_i \right) + \frac{3 c_2 k \log (n/k)}{n} + \frac{3 c_3 \xi}{n}\\
&\leq 3 (1+2c_1) \max_{1\leq i\leq k} h^2\left( Q_i, \overline{Q}_i \right) + \frac{3(1+2c_1)}{n} \sum\limits_{j=1}^n h^2\left( \sum\limits_{i=1}^k \overline{\beta}_i \overline{Q}_i, P_j \right)\\
&+ \frac{3 c_2 k \log (n/k)}{n} + \frac{3 c_3 \xi}{n},
\end{align*}
using (\ref{eq:mix_hellinger_upper}), and we can conclude with Lemma \ref{lem:aux_1}.
\section{Proofs of Section \ref{sec:main}}
\subsection{Proof of Proposition \ref{prop:mle_predictor}}
\label{sec:proof_prop_mle_predictor}
It is a direct consequence of Corollary 1 of \cite{baraudrevisited} since
\[ \overline{\mathcal{M}}(f,\alpha) := \left\{ \sum\limits_{i=1}^k \beta_i \alpha_i^{-1} f_i; \beta \in \mathcal{W}_k \right\} \]
is a convex set of densities.
\subsection{Proof of Theorem \ref{th:predictor}}
\label{sec:proof_th_predictor}
From Lemma \ref{lem:aux_2}, for all $\xi>0$ and all $\overline{\beta}\in\mathcal{W}_k$, we have
\begin{align*}
\sum\limits_{j=1}^n h^2\left( P_j, \hat{P} \right) &\leq c_1 \sum\limits_{j=1}^n h^2\left( P_j, \sum\limits_{i=1}^k \overline{\beta}_i \alpha_i^{-1} f_i \cdot\mathrm{d}\mu \right) + c_2 k\log (n/k) + c_3 \xi,
\end{align*}
with probability at least $1-e^{-\xi}$. Therefore, we have
\begin{align*}
h^2\left( \sum\limits_{i=1}^k \overline{\beta}_i \alpha_i^{-1} f_i \cdot \mathrm{d}\mu, \sum\limits_{i=1}^k \hat{\beta}_i \alpha_i^{-1} f_i \cdot \mathrm{d}\mu \right) &\leq \frac{2}{n} \sum\limits_{j=1}^n h^2\left( \sum\limits_{i=1}^k \overline{\beta}_i \alpha_i^{-1} f_i \cdot \mathrm{d}\mu, P_j \right)\\
&+ \frac{2}{n} \sum\limits_{j=1}^n h^2\left( P_j, \sum\limits_{i=1}^k \hat{\beta}_i \alpha_i^{-1} f_i \cdot \mathrm{d}\mu \right)\\
&\leq \frac{2(1+c_1)}{n} \sum\limits_{j=1}^n h^2\left( P_j, \sum\limits_{i=1}^k \overline{\beta}_i Q_i \right)\\
&+ \frac{2 c_2 k \log (n/k)}{n} + \frac{2 c_3 \xi}{n},
\end{align*}
with probability at least $1-e^{-\xi}$. We can conclude with Lemma \ref{lem:aux_1}. The second inequality can be obtained following the proof of Theorem \ref{th:conditional}.
\subsection{Proof of Proposition \ref{prop:confusion_invertible}}
\label{sec:proof_prop_confusion_invertible}
\begin{itemize}
\item If $M(f)$ is not invertible, there is $v\neq 0$ in $\mathbb{R}^k$ such that
\[ \sum\limits_{i=1}^k v_i M(f)_{ij} = \int \sum\limits_{i=1}^k v_i f_i(x) Q^*_j(dx) = 0, \]
for all $j$. Therefore, for $Q^*_j$-almost all $x$, we have $\sum\limits_{i=1}^k v_i f_i(x) =0$ for all $j$. In particular, it means that $\sum\limits_{i=1}^k v_i f_i(x) =0$, for $Q^*_{\Sigma}$-almost all $x$, i.e. the measures
\[ f_1\cdot \mathrm{d}Q^*_{\Sigma}, \dots, f_k\cdot \mathrm{d}Q^*_{\Sigma} \]
are linearly dependent.
\item Note that the distributions $Q^*_1,\dots,Q^*_k$ can be expressed as in (\ref{eq:f_to_q}). If $Q^*_1,\dots,Q^*_k$ are linearly dependent, there is $v\neq 0$ in $\mathbb{R}^k$ such that
\[ \sum\limits_{i=1}^k v_i \alpha_i^{-1} f^{\alpha}_i(x) = 0 \]
for $P_{Q^*,\alpha}$-all $x$. Since $Q^*_{\Sigma}\ll P_{Q^*,\alpha}$, we have that
\[ \sum\limits_{i=1}^k w_i f^{\alpha}_i(x) = 0 \]
for $Q^*_{\Sigma}$-all $x$, or equivalently $\sum\limits_{i=1}^k w_i f^{\alpha}_i \cdot \mathrm{d}Q^*_{\Sigma} = 0$, where $w\neq 0$ is given by $w_i=v_i/\alpha_i$ for all $i$.
\item If $Q^*_1,\dots,Q^*_k$ are linearly independent. Let $v\in\mathbb{R}^k$ be such that $0 = v_1 f^{\alpha}_1 \cdot \mathrm{d}Q^*_{\Sigma} + \dots + v_k f^{\alpha}_k \cdot \mathrm{d}Q^*_{\Sigma}$ or equivalently
\[ 0 = \sum\limits_{i=1}^k v_i f^{\alpha}_i(x), \]
for $Q^*_{\Sigma}$-almost all $x$. Since $P_{Q^*,\alpha}\ll Q^*_{\Sigma}$ we have $0 = \sum\limits_{i=1}^k v_i f^{\alpha}_i(x)$ for $P_{Q^*,\alpha}$-almost all $x$ or equivalently
\[ 0 = \sum\limits_{i=1}^k v_i f^{\alpha}_i \cdot \mathrm{d}P_{Q^*,\alpha} = \sum\limits_{i=1}^k v_i \alpha_i Q^*_i. \]
By linear independence we must have $v_i\alpha_i=0$ for all $i\in[k]$. Since $\alpha\in\mathcal{W}_k^*$, we have $v=0$ which shows that the distributions $f^{\alpha}_1\cdot \mathrm{d}Q^*_{\Sigma},\dots,f^{\alpha}_k\cdot \mathrm{d}Q^*_{\Sigma}$ are linearly independent.
\end{itemize}
\subsection{Proof of Corollary \ref{coro:predictor}}
\label{sec:proof_coro_predictor}
Under Assumption \ref{hyp:iid_mix}, for all $\eta\in\mathcal{W}_k$ we have
\[ P^*(dx) = \sum\limits_{i=1}^k \beta^*_i \eta_i^{-1} f^{\eta}_i(x) P_{\eta}(dx) = \sum\limits_{1\leq i,j\leq k} \beta^*_i \eta_i^{-1} \eta_j f^{\eta}_i(x) Q^*_j(dx). \]
Under Assumption \ref{hyp:q_sigma}, there is $\gamma\in[0,+\infty)^k$ such that $\mu=\sum\limits_{i=1}^k \gamma_i Q^*_i\in\mathscr{P}^*(f,\alpha)$. We have
\begin{align*}
h^2\left(P^*, \sum\limits_{i=1}^k \beta^*_i \alpha_i^{-1} f_i \cdot \mathrm{d}\mu \right) &\leq d_{TV}\left( \sum\limits_{1\leq i,j\leq k} \beta^*_i \alpha_i^{-1} \alpha_j f^{\alpha}_i \cdot \mathrm{d}Q^*_j, \sum\limits_{1\leq i,j\leq k} \beta^*_i \alpha_i^{-1} \gamma_j f_i \cdot \mathrm{d}Q^*_j \right)\\
&\leq \frac{1}{2} \sum\limits_{1\leq i,j\leq k} \beta^*_i \alpha_i^{-1} \alpha_j \int  |f^{\alpha}_i - f_i| Q^*_j + \frac{1}{2} \sum\limits_{1\leq i,j\leq k} \beta^*_i \alpha_i^{-1} |\alpha_j - \gamma_j|\\
&\leq \frac{\max\limits_{1\leq i\leq k} (\beta^*_i/\alpha_i) }{2} \left( \mathbb{E}_{P_{Q^*,\alpha}}\left[ ||f^{\alpha}-f||_1 \right] + ||\alpha-\gamma||_1 \right).
\end{align*}
\subsection{Proof of Corollary \ref{coro:predictor_calibration}}
\label{sec:proof_coro_predictor_calibration}
Note that Assumption \ref{hyp:q_sigma} and the fact that $f$ is calibrated implies that $P_{Q^*,\alpha}$ belongs to $\mathscr{P}^*(f,\alpha)$. We have
\begin{align*}
h^2\left(P^*, \sum\limits_{i=1}^k \beta^*_i \alpha_i^{-1} f_i \cdot \mathrm{d}P_{Q^*,\alpha} \right) &\leq d_{TV}\left( \sum\limits_{i=1}^k \beta^*_i \alpha_i^{-1} f^{\alpha}_i(x) P_{Q^*,\alpha}(dx), \sum\limits_{i=1}^k \beta^*_i \alpha_i^{-1} f_i(x) P_{Q^*,\alpha} \right)\\
&= \frac{1}{2} \displaystyle\int  \left| \sum\limits_{i=1}^k \beta^*_i \alpha_i^{-1} (f^{\alpha}_i - f_i)(x) \right| P_{Q^*,\alpha}(dx).
\end{align*}
From Lemma \ref{lem:aux_2}, with probability at least $1-e^{-\xi}$, we have
\begin{align*}
h^2\left( \sum\limits_{i=1}^k \beta^*_i \alpha_i^{-1} f_i \cdot \mathrm{d}P_{Q^*,\alpha}, \sum\limits_{i=1}^k \hat{\beta}_i \alpha_i^{-1} f_i \cdot \mathrm{d}P_{Q^*,\alpha} \right) &\leq 2 h^2\left( \sum\limits_{i=1}^k \beta^*_i \alpha_i^{-1} f_i \cdot \mathrm{d}P_{Q^*,\alpha}, P^* \right)\\
&+ 2 h^2\left( P^*, \sum\limits_{i=1}^k \hat{\beta}_i \alpha_i^{-1} f_i \cdot \mathrm{d}P_{Q^*,\alpha} \right)\\
&\leq 2(1+c_1)h^2\left( \sum\limits_{i=1}^k \beta^*_i \alpha_i^{-1} f_i \cdot \mathrm{d}P_{Q^*,\alpha}, P^* \right)\\
&+ 2c_2 \frac{k\log (n/k)}{n} + 2 c_3 \frac{\xi}{n},
\end{align*} 
for all $\xi>0$. We can conclude with Lemma \ref{lem:aux_1}.
\section{Auxiliary results}
\subsection{Proof of Lemma \ref{lem:aux_2}}
\label{sec:proof_lem_aux2}
The result combines Theorem 1 of \cite{baraudrevisited} and a bound on the $\rho$-dimension that differs slightly from the results in \cite{lecestre_mixture}. The notion of $\rho$-dimension function is introduced in Section \ref{sec:rho_dim}. From \cite{baraudrevisited}, we have
\[ \sum\limits_{j=1}^n h^2\left( P_j, \hat{P} \right) \leq a_1 \sum\limits_{j=1}^n h^2\left( P_j, \sum\limits_{i=1}^k \overline{w}_i Q_i \right) + a_2 \left( \frac{D^{\mathscr{M}_{mix}(Q_1,\dots,Q_k)}(\mathbf{P},\overline{P})}{4.7} + 1.49 + \xi \right), \]
with probability at least $1-e^{-\xi}$ for all $\xi>0$, where $a_1=150$, $a_2=5014$, $\mathbf{P}=\bigotimes_{i=1}^n P_i$, $\overline{P}\in\mathscr{M}_{mix}(Q_1,\dots,Q_k)$, and $D^{\mathscr{M}_{mix}(Q_1,\dots,Q_k)}$ is the $\rho$-dimension function associated to the model $\mathscr{M}_{mix}(Q_1,\dots,Q_k)$. The constants $a_1$ and $a_2$ are given in the proof of Theorem 1 in \cite{juntong_yannick} on page 32. The following result gives a bound on the $\rho$-dimension when we consider a class of density functions that is VC-subgraph. We refer to \cite{VanDerVaart} (Section 2.6) and \cite{rho_inventiones} (Section 8) for more on the topic of VC-classes of functions.
\begin{prop}
\label{prop:rho_dim_vc}
Let $\mathscr{F}$ be a countable subset of $\mathscr{P}_X$ and $\mathcal{F}$ an associated (countable) class of densities with respect to a $\sigma$-finite measure $\mu$, i.e. $\mathscr{F}=\{ f\cdot \mathrm{d}\mu: f\in\mathcal{F} \}$. If $\mathcal{F}$ is VC-subgraph with VC-dimension not larger than $V$, for all $\mathbf{P}=P_1\otimes\dots\otimes P_n \in\mathscr{P}_X^{\otimes n}$ and all $\overline{P}\in\mathscr{P}_X$, we have 
\[ D^{\mathscr{F}}\left(\mathbf{P},\overline{P}\right) \leq 91 \sqrt{2} V \left[ 9.11 + \log_+\left( \frac{n}{V } \right) \right]. \]
where $\log_+(x)=\max(0,\log x)$ for all $x>0$ and $D^{\mathscr{F}}$ is the $\rho$-dimension function introduced in Section \ref{sec:rho_dim}.
\end{prop}
This result is proven in Section \ref{sec:proof_rho_dim_vc}. From Lemma 2.6.15 in \cite{VanDerVaart}, the class of density functions $\mathcal{M}_{mix}(q_1,\dots,q_k)$ given by (\ref{eq:mix_fix_dens}) is VC-subgraph with VC-dimension smaller than or equal to $k+1$. Therefore we have
\[ D^{\mathscr{M}_{mix}(Q_1,\dots,Q_k)}\left( \mathbf{P}, \overline{P} \right) \leq  91 \sqrt{2} (k+1) \left[  9.11 + \log_+\left( \frac{n}{k+1} \right) \right], \]
for all $P_1,\dots,P_n,\overline{P}\in\mathscr{P}_X$. Since $k\geq 2$ and $n \geq e\times k$, we have
\[ D^{\mathscr{M}_{mix}(Q_1,\dots,Q_k)}\left( \mathbf{P}, \overline{P} \right) \leq \zeta k \log (n/k), \]
where $\zeta = 91\sqrt{2} \times \frac{3}{2} \times 10.11$. Since $a_2<5014$, we can conclude with
\[ a_2 \times \left(\frac{\zeta}{4.7}+1.49\right) < 2.1 \times 10^6. \]
\subsection{The \texorpdfstring{$\rho$}{rho}-dimension function}
\label{sec:rho_dim}
The $\rho$-dimension function is properly defined in \cite{baraudrevisited}. We slightly modify and adapt original definitions to our context in order to simplify them. One can check that the function $\psi$ defined by (\ref{eq:psi_rho}) satisfies Assumption 2 of \cite{baraudrevisited} with $a_0 = 4, a_1 = 3/8$ and $a_2^2 = 3\sqrt{2}$ (see Proposition 3 of \cite{baraudrevisited}) which gives the different constants. Let $\mathscr{M}$ be a countable subset of $\mathscr{P}_X$. For $y>0$, $P_1,\dots,P_n\in\mathscr{P}_X$ and $\overline{P}\in\mathscr{P}_X$ we write
\[ \mathscr{B}^{\mathscr{\mathscr{M}}}\left(\mathbf{P},\overline{P},y\right) := \left\{ Q\in\mathscr{M}; \sum\limits_{i=1}^n h^2\left( P_i, \overline{P} \right) + \sum\limits_{i=1}^n h^2\left( P_i, Q \right) < y^2 \right\}, \]
where $\mathbf{P}$ is the product distribution $P_1\otimes\dots\otimes P_n$. If $\mathcal{M}$ is a set of probability density functions with respect to a $\sigma$-finite measure $\mu$ such that
\begin{equation}
\label{eq:representation}
\mathscr{M}\cup\{\overline{P}\}=\{ q\cdot \mathrm{d}\mu; q\in\mathcal{M}\},
\end{equation}
we write
\[ w\left(\mu,\mathcal{M},\mathscr{M},\mathbf{P},\overline{P},y\right) := \left[ \sup\limits_{Q\in\mathscr{B}^{\mathscr{M}}(\overline{P},y)} \left| \mathbf{T}\left(\mathbf{X},\overline{p},q\right) - \mathbb{E}_{\mathbf{P}} \left[ \mathbf{T}\left(\mathbf{X},\overline{p},q\right) \right] \right| \right]. \]
Similarly, we define
\[ \mathbf{w}^{\mathscr{M}}\left(\mathbf{P},\overline{P},y\right) = \inf_{(\mu,\mathcal{M})} w\left(\mu,\mathcal{M},\mathscr{M},\mathbf{P},\overline{P},y\right), \]
where the infimum is taken over all couples $(\mu,\mathcal{M})$ satisfying (\ref{eq:representation}). We can now define the $\rho$-dimension function $D^{\mathscr{M}}$ by 
\[ D^{\mathscr{M}}\left(\mathbf{P},\overline{P} \right) := \left[ \upsilon \sup\left\{ y^2 ;\mathbf{w}^{\mathscr{M}}\left(\overline{P},y\right)>\omega y^2 \right\} \right] \bigvee 1, \]
with $\upsilon=3/2^{10+1/2}$ and $\omega=3/64$.
\subsubsection{Proof of Proposition \ref{prop:rho_dim_vc}}
\label{sec:proof_rho_dim_vc}
Since $\mathcal{F}$ is VC-subgraph, the set $\left\{ \psi\left( \sqrt{\frac{f}{\overline{p}}} \right); f\in\mathcal{F} \right\}$ is also VC-subgraph with VC-dimension not larger than $V$ (see proof of Proposition 42 (vii) in \cite{rho_inventiones}), such as any of its subsets. In particular, we can consider
\[ \mathcal{F}(\mathbf{P},\overline{P},y) = \left\{ \psi\left( \sqrt{\frac{f}{\overline{p}}} \right); f\in\mathcal{F}, \sum\limits_{j=1}^n h^2\left(P_j,\overline{P}\right) + h^2\left(P_j,F\right) < y^2 \right\}. \]
From Theorem 2 in \cite{juntong_yannick} and (11) in \cite{baraudrevisited}, we have
\[ \mathbb{E}\left[ \sup_{f\in\mathcal{F}(\mathbf{P},\overline{P},y)} \left|\sum\limits_{j=1}^n \psi(\sqrt{f/\overline{p}}) - \mathbb{E}\left[ \psi(\sqrt{f/\overline{p}}) \right] \right| \right] \leq 4.74 \sqrt{V y^2 a_2^2 \mathcal{L}(y)} + 90 V \mathcal{L}(y), \]
where $\mathcal{L}(y) = 9.11+\log_+( n / y^2 a_2^2)$.
We can now follow the structure of the proof of Proposition 7 in \cite{juntong_yannick}. With the notation of \cite{baraudrevisited}, we have
\begin{align*}
w^{\mathscr{F}}\left( \mathbf{P},\overline{P}^{\otimes n},y \right) &= \mathbb{E}\left[ \sup_{f\in\mathcal{F}(\mathbf{P},\overline{P},y)} \left|\sum\limits_{j=1}^n \psi(\sqrt{f/\overline{p}}) - \mathbb{E}\left[ \psi(\sqrt{f/\overline{p}}) \right] \right| \right]\\
&\leq 4.74 a_2 y \sqrt{V \mathcal{L}(y)} + 90 V \mathcal{L}(y).
\end{align*}
Let $D\geq a_1^2 V /(16 a_2^4) = 2^{-11} V$ to be chosen later on and $\beta = a_1/(4a_2)$. For $y\geq \beta^{-1} \sqrt{D}$,
\[ L(y) = 9.11 + \log_+\left( \frac{n}{y^2 a_2^2} \right) \leq 9.11 +\log_+\left( \frac{n}{V} \right) = L. \]
Hence for all $y\geq \beta^{-1}\sqrt{D}$,
\begin{align*}
w^{\mathscr{F}}\left( \mathbf{P},\overline{P}^{\otimes n},y \right) &\leq 4.74 a_2 y \sqrt{V L} + 90 V L\\
&\leq \frac{a_1 y^2}{8} \left[ 1.185 \frac{\sqrt{V L}}{\sqrt{D}} + \frac{45}{\sqrt{2}} \frac{V L}{D} \right]\\
&\leq \frac{a_1 y^2}{8},
\end{align*}
for $D = 91 \sqrt{2} VL > 2^{-11} V$. The result follows from the definition of the $\rho$-dimension given in \cite{baraudrevisited} (Definition
4).
\section{Proofs of Section \ref{sec:calibration}}
\subsection{Proof of Lemma \ref{lem:canonical_calibration_exp}}
\label{sec:proof_lem_canonical_calibration_exp}
Since $f$ is canonically calibrated with respect to $\Pi$, we have
\begin{align*}
\mathbb{E}_{P_{Q,\beta}}\left[ \phi(f(X)) \right] &= \sum\limits_{i\in I} \beta_i \mathbb{E}_{Q_i}\left[ \phi(f(X)) \right]\\
&= \sum\limits_{i\in I} \beta_i \alpha_i^{-1} \mathbb{E}_{\Pi}\left[ \mathbbm{1}_{Y=i} \phi(f(X)) \right]\\
&= \sum\limits_{i\in I} \beta_i \alpha_i^{-1} \mathbb{E}_{\Pi}\left[ f_i(X) \phi(f(X)) \right]\\
&= \sum\limits_{i\in I} \beta_i \alpha_i^{-1} \mathbb{E}_{P_{Q,\alpha}}\left[ f_i(X) \phi(f(X)) \right]\\
&= \mathbb{E}_{P_{Q,\alpha}}\left[ \sum\limits_{i=1}^k \beta_i \alpha_i^{-1}  f_i(X) \phi(f(X)) \right]\\
&= \mathbb{E}_{P_{f,\beta}}\left[ \phi(f(x)) \right],
\end{align*}
\subsection{Proof of Lemma \ref{lem:calibration_product}}
\label{sec:proof_lem_calibration_product}
Canonical calibration means
\begin{align*}
& g_{i,j}(X_1,X_2) = \mathbb{E}_{\Pi}\left[ \mathbbm{1}_{(Y_1,Y_2)=(i,j)} | g(X_1,X_2) \right]\\
\Leftrightarrow & f_i(X_1) f_j(X_2) = \mathbb{E}_{\Pi}\left[ \mathbbm{1}_{Y_1=i} \mathbbm{1}_{Y_2=j} | f^{(1)}(X_1), f^{(2)}(X_2) \right]
\end{align*}
for all $i\in I_1$ and $j\in I_2$. Since $f_1$ and $f_2$ are canonically calibrated, it suffices to show that
\[ \mathbb{E}_{\Pi_1}\left[ \mathbbm{1}_{Y_1=i} | f^{(1)}(X_1) \right] \mathbb{E}_{\Pi_2}\left[ \mathbbm{1}_{Y_2=j} | f^{(2)}(X_2) \right] = \mathbb{E}_{\Pi}\left[ \mathbbm{1}_{(Y_1,Y_2)=(i,j)} | f^{(1)}(X_1), f^{(2)}(X_2) \right] \]
We can do that with standard results from measure theory. Since $(X_1,Y_1)$ and $(X_2,Y_2)$ are independent, we have
\begin{align*}
\mathbb{E}_{\Pi}\left[ \mathbbm{1}_{(Y_1,Y_2)=(i,j)} | f^{(1)}(X_1), f^{(2)}(X_2) \right] &= \mathbb{E}_{\Pi}\left[ \mathbbm{1}_{Y_1=i} | f^{(1)}(X_1),f^{(2)}(X_2) \right] \mathbb{E}_{\Pi}\left[ \mathbbm{1}_{Y_2=j} | f^{(1)}(X_1), f^{(2)}(X_2) \right]\\
&= \mathbb{E}_{\Pi}\left[ \mathbbm{1}_{Y_1=i} | f^{(1)}(X_1) \right] \mathbb{E}_{\Pi}\left[ \mathbbm{1}_{Y_2=j} | f^{(2)}(X_2) \right].
\end{align*}
\subsection{Proof of Proposition \ref{prop:kl_min_canon}}
\label{sec:proof_prop_kl_min_canon}
With Lemma \ref{lem:canonical_calibration_exp}, we have
\begin{align*}
\mathbb{E}_{P^*}\left[ \log\left( \sum\limits_{i=1}^k \beta_i \alpha_i^{-1} f_i(x) \right) \right] &= \mathbb{E}_{P_{\beta^*}}\left[ \log\left( \sum\limits_{i=1}^k \beta_i \alpha_i^{-1} f_i(x) \right) \right]\\
&= \mathbb{E}_{P_{f,\beta^*}}\left[ \log\left( \sum\limits_{i=1}^k \beta_i \alpha_i^{-1} f_i(x) \right) \right]\\
&= \int \log\left( \sum\limits_{i=1}^k \beta^*_i \alpha_i^{-1} f_i(x) \right) P_{f,\beta^*} - \mathbf{K}(P_{f,\beta^*}||P_{f,\beta}).
\end{align*}
We can conclude with the fact that the first term on the right hand side does not depend on $\beta$.
\subsection{Proof of Theorem \ref{th:calibration}}
\label{sec:proof_th_calibration}
We break down the proof into several steps. We first get a deviation inequality for the Hellinger distance between $P_{f,\alpha,\beta}$ and $P_{f,\alpha,\hat{\beta}}$ when we assume that the $X_1,\dots,X_n$ are i.i.d. with distribution $P_{f,\alpha,\beta}$. Then we use the fact that $f$ is canonically calibrated to get rid of the deviation term.
\begin{itemize}
    \item Lemma \ref{lem:marginal_calibration} implies that $P_{Q^*,\alpha}\in\mathscr{P}(f,\alpha)$ since $f$ is canonically calibrated with respect to $\Pi_{\alpha}$. We consider the $\rho$-estimator $\hat{P}=\hat{P}(\mathbf{X},\mathcal{M}(f,\alpha))\in \overline{M}(f,\alpha,P_{Q^*,\alpha})$ defined in Section \ref{sec:rho_est_pred} and following (\ref{eq:hat_beta_predictor}) we denote by $\hat{\beta}$ any element of $\mathcal{W}_k$ such that
\[ \hat{P} = \sum\limits_{i=1}^k \hat{\beta}_i \alpha_i^{-1} f_i \cdot \mathrm{d}P_{Q^*,\alpha}=P_{f,\alpha,\hat{\beta}}. \]
    \item Let $\beta$ be in $\mathcal{W}_k$. If we apply Lemma \ref{lem:aux_2} with $P_1=\dots=P_n=P_{f,\alpha,\beta}$ and $\overline{w}=\beta$, we have
\begin{equation}
\label{eq:calibration_proof_1}
\mathbb{P}_{\mathbf{X}\sim P_{f,\alpha,\beta}^{\otimes n}}\left( h^2\left( P_{f,\alpha,\beta}, P_{f,\alpha,\hat{\beta}} \right) \leq 0 + c_2 \frac{k\log(n/k)}{n} + c_3 \frac{\xi}{n} \right) \geq 1-e^{-\xi},
\end{equation}
for all $\xi>0$. This what we get if we assume that the $X_1,\dots,X_n$ are i.i.d. with distribution $P_{f,\alpha,\beta}$. 
    \item Corollary \ref{coro:calibration_product} implies that the function
    \[ g: \left| \begin{array}{ll}
		\mathscr{X}^n &\rightarrow \mathcal{W}_{I^n}\\
        (x_1,\dots,x_n) &\mapsto (f(x_1)_{i_1}\dots f(x_n)_{i_n})_{i\in I^n}
	\end{array} \right. \]
    is canonically calibrated with respect to the distribution $\Pi^{(n)}$ defined by
    \[ (X_1,Y_1\dots,X_n,Y_n) \sim \Pi_{\alpha}^{\otimes n} \Leftrightarrow (X_1,\dots,X_n,Y_1,\dots,Y_n) \sim \Pi^{(n)}. \]
    \item One can check from the definition and the chosen model that $\hat{P}$ is a measurable function of $(f(X_1),\dots,f(X_n))$, and equivalently it is a function of $g(\mathbf{X})$. We now use Assumption \ref{hyp:iid_mix},  Using Lemma \ref{lem:canonical_calibration_exp} and (\ref{eq:calibration_proof_1}), we have
\begin{align*}
&\mathbb{P}_{\mathbf{X}\sim (P^*)^{\otimes n}}\left( h^2\left( P_{f,\alpha,\beta}, P_{f,\alpha,\hat{\beta}} \right) \leq c_2 \frac{k\log(n/k)}{n} + c_3 \frac{\xi}{n} \right)\\
&= \mathbb{P}_{(X_1,Y_1)\dots(X_n,Y_n)\sim \Pi^{(n)}}\left( h^2\left( P_{f,\alpha,\beta}, P_{f,\alpha,\hat{\beta}} \right) \leq c_2 \frac{k\log(n/k)}{n} + c_3 \frac{\xi}{n} \right)\\
&= \mathbb{P}_{(X_1,Y_1)\dots(X_n,Y_n)\sim \Pi^{(n)}}\left( h^2\left( P_{f,\alpha,\beta}, P_{f,\alpha,\hat{\beta}} \right) \leq c_2 \frac{k\log(n/k)}{n} + c_3 \frac{\xi}{n} \right)\\
&= \mathbb{P}_{\mathbf{X}\sim P_{f,\alpha,\beta}^{\otimes n}}\left( h^2\left( P_{f,\alpha,\beta}, P_{f,\alpha,\hat{\beta}} \right) \leq c_2 \frac{k\log(n/k)}{n} + c_3 \frac{\xi}{n} \right) \geq 1-e^{-\xi}, 
\end{align*}
for all $\xi>0$. This proves (\ref{eq:th_calibration_1}) and (\ref{eq:th_calibration_2}) is obtained with Lemma \ref{lem:aux_1} and the fact that $P_{Q^*,\alpha}\in\mathscr{P}^*(f,\alpha)$.
\end{itemize}
\end{appendix}

\end{document}